\documentclass[english,reqno]{amsart}

\usepackage[margin=2.5cm]{geometry}

\input{preamble}
\newcommand{\RNum}[1]{\uppercase\expandafter{\romannumeral #1\relax}}

\newcommand{\Cbb}{\mathbb{C}}
\newcommand{\Ebb}{\mathbb{E}}
\newcommand{\Nbb}{\mathbb{N}}

\newcommand{\Rbb}{\mathbb{R}}

\newcommand{\Zbb}{\mathbb{Z}}

\newcommand{\Abf}{\mathbf{A}}
\newcommand{\bbf}{\mathbf{b}}

\newcommand{\Gbf}{\mathbf{G}}

\newcommand{\Hbf}{\mathbf{H}}

\newcommand{\Ibf}{\mathbf{I}}

\newcommand{\Lbf}{\mathbf{L}}
\newcommand{\mbf}{\mathbf{m}}
\newcommand{\Mbf}{\mathbf{M}}
\newcommand{\nbf}{\mathbf{n}}

\newcommand{\pbf}{\mathbf{p}}

\newcommand{\sbf}{\mathbf{s}}

\newcommand{\tbf}{\mathbf{t}}

\newcommand{\ubf}{\mathbf{u}}
\newcommand{\Ubf}{\mathbf{U}}

\newcommand{\Vbf}{\mathbf{V}}
\newcommand{\wbf}{\mathbf{w}}
\newcommand{\Wbf}{\mathbf{W}}
\newcommand{\xbf}{\mathbf{x}}
\newcommand{\Xbf}{\mathbf{X}}
\newcommand{\ybf}{\mathbf{y}}
\newcommand{\Ybf}{\mathbf{Y}}
\newcommand{\zbf}{\mathbf{z}}

\newcommand{\ep}{\epsilon}

\newcommand{\CalA}{{\mathcal{A}}}

\newcommand{\CalF}{{\mathcal{F}}}

\newcommand{\CalI}{{\mathcal{I}}}

\newcommand{\CalL}{{\mathcal{L}}}
\newcommand{\CalN}{{\mathcal{N}}}
\newcommand{\CalO}{{\mathcal{O}}}
\newcommand{\CalP}{{\mathcal{P}}}

\newcommand{\CalS}{{\mathcal{S}}}

\newcommand{\nrm}[1]{\left\Vert {#1} \right\Vert}

\newcommand{\supp}[1]{\text{supp}\left(#1\right)}

\newcommand{\txt}[2]{\hspace{#1cm} \text{#2} \hspace{#1cm}} 
\newcommand{\symb}[2]{\hspace{#1cm} {#2} \hspace{#1cm}}
\renewcommand{\tilde}[1]{\widetilde{#1}}

\newcommand{\argmin}{\text{argmin}}

\newcommand{\lan}{\left\langle}
\newcommand{\ran}{\right\rangle}

\newcommand{\wstar}{\wbf^\star}
\newcommand{\what}{{\widehat{\wbf}}}

\usepackage{array}
\begin{document}

%%%%%%%%%%%%%%%%%%% Publisher's Area please ignore %%%%%%%%%%%%%%%%%%%%%%%
%
%\catchline{}{}{}{}{}
%
%%%%%%%%%%%%%%%%%%%%%%%%%%%%%%%%%%%%%%%%%%%%%%%%%%%%%%%%%%%%%%%%%%%%%%%%%%

\title{Weak SINDy for Partial Differential Equations}

\author{Daniel A. Messenger, David M. Bortz}
\email{daniel.messenger@colorado.edu, dmbortz@colorado.edu\footnote{Department of Applied Mathematics, University of Colorado Boulder, 11 Engineering Dr., Boulder, CO 80309, USA.}}
\maketitle

%\begin{history}
%\received{(Day Month Year)}
%\revised{(Day Month Year)}
%\accepted{(Day Month Year)}
%\comby{(xxxxxxxxxx)}
%\end{history}

\begin{abstract}
Sparse Identification of Nonlinear Dynamics (SINDy) is a method of system discovery that has been shown to successfully recover governing dynamical systems from data \cite{brunton2016discovering,rudy2017data}. Recently, several groups have independently discovered that the weak formulation provides orders of magnitude better robustness to noise. Here we extend our Weak SINDy (WSINDy) framework introduced in \cite{messenger2020weak} to the setting of partial differential equations (PDEs). The elimination of pointwise derivative approximations via the weak form enables effective machine-precision recovery of model coefficients from noise-free data (i.e.\ below the tolerance of the simulation scheme) as well as robust identification of PDEs in the large noise regime (with signal-to-noise ratio approaching one in many well-known cases). This is accomplished by discretizing a convolutional weak form of the PDE and exploiting separability of test functions for efficient model identification using the Fast Fourier Transform. The resulting WSINDy algorithm for PDEs has a worst-case computational complexity of $\CalO(N^{D+1}\log(N))$ for datasets with $N$ points in each of $D+1$ dimensions (i.e.\ $\CalO(\log(N))$ operations per datapoint). Furthermore, our Fourier-based implementation reveals a connection between robustness to noise and the spectra of test functions, which we utilize in an \textit{a priori} selection algorithm for test functions. Finally, we introduce a learning algorithm for the threshold in sequential-thresholding least-squares (STLS) that enables model identification from large libraries, and we utilize scale-invariance at the continuum level to identify PDEs from poorly-scaled datasets. We demonstrate WSINDy's robustness, speed and accuracy on several challenging PDEs.
\end{abstract}


\section{Introduction}\label{sec:intro}

Stemming from Akaike's seminal work in the 1970's \cite{Akaike1974IEEETransAutomControl,Akaike1977Applicationsofstatistics}, research into the automatic creation of accurate mathematical models from data has progressed dramatically. In the last 20 years, substantial developments have been made at the interface of applied mathematics and statistics to design data-driven model selection algorithms that are both statistically rigorous and computationally efficient (see \cite{BortzNelson2006BullMathBiol,LagergrenNardiniMichaelLavigneEtAl2020ProcRSocA,LillacciKhammash2010PLoSComputBiol,ToniWelchStrelkowaEtAl2009JRSocInterface,WarneBakerSimpson2019BullMathBiola,WuWu2002StatistMed}
for both theory and applications). An important achievement in this field was the formulation and subsequent discretization of the system discovery problem in terms of a candidate basis of nonlinear functions evaluated at the given dataset, together with a sparsification measure to avoid overfitting \cite{crutchfield1987equations}. In \cite{wang2011predicting} the authors extended this framework to the context of catastrophe prediction and used compressed sensing techniques to enforce sparsity. More recently, this approach has been generalized as the SINDy algorithm (Sparse Identification of Nonlinear Dynamics) \cite{brunton2016discovering} and successfully used to identify a variety of discrete and continuous dynamical systems.

The wide applicability, computational efficiency, and interpretability of the SINDy algorithm has spurred an explosion of interest in the problem of identifying nonlinear dynamical systems from data \cite{cortiella2020sparse,qin2020data,dai2020detecting,de2020pysindy,hoffmann2019reactive,tran2017exact,mangan2017model}. In addition to the sparse regression approach adopted in SINDy, some of the primary techniques include Gaussian process regression \cite{owhadi2015bayesian,raissi2017machine}, deep neural networks \cite{rudy2019deep,wang2020deep,LagergrenNardiniMichaelLavigneEtAl2020ProcRSocA}, Bayesian inference \cite{zhang2018robust,zhang2019robust,wang2020perspective} and classical methods from numerical analysis \cite{kang2019ident,keller2019discovery,wu2019numerical}. The variety of approaches for model discovery from data qualitatively differ in the interpretability of the resulting data-driven dynamical system, the computational efficiency of the algorithm, and the robustness to noise, scale separation, etc. For instance, a neural-network based data-driven dynamical system does not easily lend itself to physical interpretation\footnote{There have been efforts to address the interpretability of neural networks, see e.g.
\cite{montavon2018methods,toms2020physically,rudin2019stop}.}. The SINDy algorithm allows for direct interpretations of the dynamics from identified differential equations and uses sequential-thresholding least-squares (STLS) to enforce a sparse solution $x\in \Rbb^n$ to a linear system $Ax=b$. STLS has been proven to converge to a local minimizer of the non-convex functional $F(x) = \nrm{Ax-b}_2^2+\lambda^2\nrm{x}_0$ in at-most $n$ iterations \cite{zhang2019convergence}. 

The aim of the present article is to extend the Weak SINDy method (WSINDy) for recovering ordinary differential equations (ODEs) from data to the context of partial differential equations (PDEs) \cite{messenger2020weak}. WSINDy is a Galerkin-based data-driven model selection algorithm that utilizes the weak form of the dynamics in a sparse regression framework. By integrating in time against compactly-supported test functions, WSINDy avoids approximation of pointwise derivatives which are known to result in low robustness to noise \cite{rudy2017data}. In \cite{messenger2020weak} we showed that by integrating against a suitable choice of test functions, correct ODE model terms can be identified together with machine-precision recovery of coefficients (i.e.\ below the tolerance of the data simulation scheme) from noise-free synthetic data, and for datasets with large noise, WSINDy successfully recovers the correct model terms without explicit data denoising. The use of integral equations for system identification was proposed as early as the 1980's \cite{crutchfield1987equations} and was carried out in a sparse regression framework in  \cite{schaeffer2017sparse} in the context of ODEs, however neither works utilized the full generality of the weak form.

Sparse regression approaches for learning PDEs from data have seen a tremendous spike in activity in the years since 2016, stemming from the pioneering works \cite{schaeffer2017learning} and \cite{rudy2017data}. The Douglas-Rachford algorithm was used in \cite{schaeffer2017learning} to enforce sparsity while \cite{rudy2017data} introduces PDE-FIND, an extension of SINDy to PDEs. Many other predominant approaches for learning dynamical systems (Gaussian processes, deep learning, Bayesian inference, etc.) have since been extended to the discovery of PDEs \cite{chen2019learning,nardini2020learning,lagergren2020biologically,LagergrenNardiniMichaelLavigneEtAl2020ProcRSocA,wang2020identification,wang2019variational,xu2020dlga,wu2020data,thaler2019sparse,WarneBakerSimpson2019BullMathBiola}.
A significant disadvantage for the vast majority of PDE discovery methods is the requirement of pointwise derivative approximations. Steps to alleviate this are taken by the authors of \cite{qin2019data} and \cite{wu2020data}, where neural network-based recovery schemes are combined with integral and abstract evolution equations to recover PDEs, and in \cite{wang2019variational}, where the finite element-based method Variational System Identification (VSI) is introduced to identify reaction-diffusion systems and uses backward Euler to approximate the time derivative.

WSINDy falls into a class of methods for discovering PDEs without any pointwise derivative approximations, black-box routines or conventional noise filtering. Through integration by parts in both space and time against smooth compactly-supported test functions, WSINDy is able to recover PDEs from datasets with much higher noise levels, and from truly weak solutions (see Figure \ref{burgerssols} in Section \ref{sec:examples}). This works suprisingly well even as the signal-to-noise ratio approaches one. Furthermore, as in the ODE setting, WSINDy achieves high-accuracy recovery in the low-noise regime. These overwhelming improvements resulting from a fully weak\footnote{The underlying true solution need only have bounded variation and the only derivatives approximated are weak derivatives.} identification method have also been discovered independently by other groups \cite{reinbold2020using,gurevich2019robust}. WSINDy offers several advantages over these alternative frameworks. Firstly, we use a convolutional weak form which enables efficient model identification using the Fast Fourier Transform (FFT). For measurement data with $N$ points in each of the $D+1$ space-time dimensions ($N^{D+1}$ total data points), the resulting algorithmic complexity of WSINDy in the PDE setting is at worst $\CalO(N^{D+1}\log(N))$, in other words $\CalO(\log(N))$ floating point operations per data-point. Subsampling further reduces the cost. Furthermore, our FFT-based approach reveals a key mechanism behind the observed robustness to noise, namely that spectral decay properties of test functions can be tuned to damp noise-dominated modes in the data, and we develop a learning algorithm for test function hyperparameters based on this mechanism. WSINDy also utilizes scale-invariance of the PDE and a modified STLS algorithm with automatic threshold selection to recover models from (i) poorly-scaled data and (ii) large candidate model libraries.

The outline of the article is as follows. In Section \ref{sec:probstatement} we define the system discovery problem that we aim to solve and the notation to be used throughout. We then introduce the convolutional weak formulation along with our FFT-based discretization in Section \ref{sec:weakdiscret}. Key ingredients of the WSINDy algorithm for PDEs (Algorithm \ref{alg:wsindypde}) are covered in Section \ref{sec:alg}, including a discussion of spectral properties of test functions and robustness to noise (\ref{sec:piecewisepoly}), our modified sequential thresholding scheme (\ref{sec:sparse}), and regularization using scale invariance of the underlying PDE (\ref{sec:rescale}). Section \ref{sec:examples} contains numerical model discovery results for a range of nonlinear PDEs, including several vast improvements on existing results in the literature. We conclude the main text in Section \ref{sec:conclusion} which summarizes the exposition and includes natural next directions for this line of research. Lastly, additional numerical details are included in the Appendix.
\section{Problem Statement and Notation}\label{sec:probstatement}
Let $\Ubf$ be a spatiotemporal dataset given on the spatial grid $\Xbf \subset \overline{\Omega}$ over timepoints $\tbf\subset [0, T]$ where $\Omega$ is an open, bounded subset in  $\Rbb^D$, $D\geq 1$. In the cases we consider here, $\Omega$ is rectangular and the spatial grid is given by a tensor product of one-dimensional grids $\Xbf = \Xbf_1\otimes \cdots \otimes \Xbf_D$, where each $\Xbf_d\in \Rbb^{N_d}$ for $1\leq d\leq D$ has equal spacing $\Delta x$, and the time grid $\tbf\in\Rbb^{N_{D+1}}$ has equal spacing $\Delta t$. The dataset $\Ubf$ is then a $(D+1)$-dimensional array with dimensions $N_1\times \cdots \times N_{D+1}$. We write $h(\Xbf,\tbf)$ to denote the $(D+1)$-dimensional array obtained by evaluating the function $h:\Rbb^D\times\Rbb\to \Cbb$ at each of the points in the computational grid $(\Xbf,\tbf)$. Individual points in $(\Xbf,\tbf)$ will often be denoted by $(\xbf_k,t_k)\in(\Xbf,\tbf)$ where
\[(\xbf_k,t_k) = (\Xbf_{k_1,\dots,k_{D}},t_{k_{D+1}}) = (x_{k_1},\dots,x_{k_D},t_{k_{D+1}})\in \Rbb^D\times \Rbb.\]
In a mild abuse of notation, for a collection of points $\{(\xbf_k,t_k)\}_{k\in[K]}\subset (\Xbf,\tbf)$, the index $k$ plays a double role as a single index in the range $[K]:=\{1,\dots,K\}$ referencing the point $(\xbf_k,t_k) \in\{(\xbf_k,t_k)\}_{k\in[K]}$ and as a multi-index on $(\xbf_k, t_k)= (\Xbf_{k_1,\dots,k_{D}}, t_{k_{D+1}})$, where $k_d$ references the $d$th coordinate. This is particularly useful for defining a matrix $\Gbf\in \Cbb^{K\times J}$ of the form 
\[\Gbf_{k,j} = h_j(\xbf_k,t_k)\]
(as in equation \eqref{conv_disc} below) where $(h_j)_{j\in[J]}$ is a collection of $J$ functions $h_j:\Rbb^D\times\Rbb\to \Cbb$ evaluated at the set of $K$ points $\{(\xbf_k,t_k)\}_{k\in[K]}\subset (\Xbf,\tbf)$. \\
 
We assume that the data satisfies $\Ubf = u(\Xbf, \tbf) + \ep$ for i.i.d.\ noise\footnote{Here  $\ep$ is used to denote a multi-dimensional array of i.i.d.\ random variables and has the same dimensions as $\Ubf$.} $\ep$ and weak solution $u$ of the PDE
\begin{equation}\label{gen_pde}
D^{\alpha^0}u(x,t) = D^{\alpha^1}g_1(u(x,t))+D^{\alpha^2}g_2(u(x,t))+\dots+D^{\alpha^S} g_S(u(x,t)), \qquad x\in \Omega,\ t\in (0,T). 
\end{equation}
The problem we aim to solve is the identification of functions $(g_s)_{s\in [S]}$ and corresponding differential operators $(D^{\alpha^s})_{s\in[S]}$ that govern the evolution\footnote{Commonly $D^{\alpha^0}$ is a time derivative $\partial_t$ or $\partial_{tt}$, although this is not required.} of $u$ according to $D^{\alpha^0}u$  given the dataset $\Ubf$ and computational grid $(\Xbf,\tbf)$. Here and throughout we use the multi-index notation $\alpha^s =~ (\alpha^s_1,\dots,\alpha^s_D,\alpha^s_{D+1}) \in \Nbb^{D+1}$ to denote partial differentiation\footnote{We will avoid using subscript notation such as $u_x$ to denote partial derivatives, instead using $D^\alpha u$ or $\partial_x u$. For functions $f(x)$ of one variable, $f^{(n)}(x)$ denotes the $n$th derivative of $f$.} with respect to $x = (x_1,\dots, x_D)$ and $t$, so that 
\[D^{\alpha^s}u(x,t) = \frac{\partial^{\alpha^s_1+\cdots+\alpha^s_D+\alpha^s_{D+1}}}{\partial x_1^{\alpha^s_1}\dots \partial x_D^{\alpha^s_D}\partial t^{\alpha^s_{D+1}}}u(x,t).\]
We emphasize that a wide variety of PDEs can be written in the form \eqref{gen_pde}. In particular, in this paper we demonstrate our method of system identification on inviscid Burgers, Korteweg-de Vries, Kuramoto-Sivashinsky, nonlinear Schr\"odinger's, Sine-Gordon, a reaction-diffusion system and Navier-Stokes. The list of admissable PDEs that can be transformed into a weak form without any derivatives on the state variables includes many other well-known PDEs (Allen-Cahn, Cahn-Hilliard, Boussinesq,\dots).
%%%
%%%
\section{Weak Formulation}\label{sec:weakdiscret}
%%%
%%%
To arrive at a computatonally tractable model recovery problem, we assume that the set of multi-indices $(\alpha^s)_{s\in [S]}$ together with $\alpha^0$ enumerates the set of possible true differential operators that govern the evolution of $u$ and that $(g_s)_{s\in [S]} \subset \text{span}(f_j)_{j\in [J]}$ where the family of functions $(f_j)_{j\in [J]}$ (referred to as the {\it trial functions}) is known beforehand. This enables us to rewrite \eqref{gen_pde} as
\begin{equation}\label{diffform}
D^{\alpha^0} u  = \sum_{s=1}^S\sum_{j=1}^J \wstar_{(s-1)J+j} D^{\alpha^s} f_j(u),
\end{equation}
so that discovery of the correct PDE is reduced to a finite-dimensional problem of recovering the vector of coefficients $\wstar\in \Rbb^{SJ}$, which is assumed to be sparse. 

To convert the PDE into its weak form, we multiply equation \eqref{diffform} by a smooth \textit{test function} $\psi(x,t)$, compactly-supported in $\Omega\times (0,T)$, and integrate over the spacetime domain,
\[\lan \psi,\ D^{\alpha^0} u\ran  = \sum_{s=1}^S\sum_{j=1}^J \wstar_{(s-1)J+j} \lan \psi,\ D^{\alpha^s} f_j(u)\ran,\]
where the $L^2$-inner product is defined $\lan \psi, f\ran := \int_0^T\int_\Omega \psi^*(x,t)f(x,t)\,dxdt$ and $\psi^*$ denotes the complex conjugate of $\psi$, although in what follows we integrate against only real-valued test functions and will omit the complex conjugation. Using the compact support of $\psi$ and Fubini's theorem, we then integrate by parts as many times as necessary to arrive at the following weak form of the dynamics:
\begin{equation}\label{weakform}
\lan (-1)^{|\alpha^0|} D^{\alpha^0}\psi,\ u \ran = \sum_{s=1}^S\sum_{j=1}^J \wstar_{(s-1)J+j} \lan(-1)^{|\alpha^s|} D^{\alpha^s}\psi,\ f_j(u)\ran,
\end{equation}
where $|\alpha^s| := \sum_{d=1}^{D+1}\alpha^s_d$ is the order of the multi-index\footnote{For example, with $D^{\alpha^s} = \frac{\partial^{2+1}}{\partial x^2\partial y}$, integration by parts occurs twice with respect to the $x$-coordinate and once with respect to $y$, so that $|\alpha^s| = 3$ and $(-1)^{|\alpha^s|}=-1$.}. Using an ensemble of test functions $(\psi_k)_{k\in [K]}$, we then discretize the integrals in \eqref{weakform} with $f_j(u)$ replaced by $f_j(\Ubf)$ (i.e. evaluated at the observed data $\Ubf$) to arrive at the linear system 
\[\bbf = \Gbf \wstar\]
defined by 
\begin{equation}\label{weak_disc}
\begin{dcases} \hspace{1.4cm}\bbf_k = \lan (-1)^{|\alpha^0|}D^{\alpha^0} \psi_k,\  \Ubf\ran,\\
\Gbf_{k,(s-1)J+j} = \lan (-1)^{|\alpha^s|}D^{\alpha^s} \psi_k,\  f_j(\Ubf)\ran,\end{dcases}
\end{equation}
where $\bbf\in\Rbb^K$, $\Gbf\in\Rbb^{K\times SJ}$ and $\wstar\in\Rbb^{SJ}$ are referred to throughout as the {\it left-hand side}, {\it Gram matrix} and {\it model coefficients}, respectively. In a mild abuse of notation, we use the inner product both in the sense of a continuous and exact integral in \eqref{weakform} and a numerical approximation in \eqref{weak_disc} which depends on a chosen quadrature rule. Building off of its success in the ODE setting, we use the trapezoidal rule throughout, as it has been shown to yield nearly negligible quadrature error with the test functions employed below (see Section \ref{sec:piecewisepoly} and \cite{messenger2020weak}). In this way, solving $\bbf = \Gbf\wstar$ for the model coefficients $\wstar$ allows for recovery of the PDE \eqref{diffform} without pointwise derivative approximations. The Gram matrix $\Gbf\in \Rbb^{K\times SJ}$ and left-hand side $\bbf\in\Rbb^K$ defined in \eqref{weak_disc} conveniently take the same form regardless of the spatial dimension $D$, as their dimensions only depend on the number of test functions $K$ and the size $SJ$ of the model library, composed of $J$ trial functions $(f_j)_{j\in[J]}$ and $S$ candidate differential operators enumerated by the multi-index set $\pmb{\alpha} := (\alpha^s)_{1\leq s\leq S}$.

%To construct and solve $\bbf = \Gbf\wstar$, we must choose the test functions $(\psi_k)_{k\in[K]}$, trial functions $(f_j)_{j\in[J]}$, multi-indices $\pmb{\alpha} := (\alpha^s)_{0\leq s\leq S}$, the quadrature rule for discretizing integrals, and the method of enforcing sparsity in the weight vector $\wstar$. The WSINDy\_PDE algorithm uses a convolution-based approach for fast computation of $\Gbf$ and $\bbf$ that operates efficiently over multi-dimensional arrays by exploiting separability in the test functions. Convolutions are computed using the Fast Fourier Transform (FFT) which also naturally relates to the selection test functions with fast-decaying spectrum over noise-dominated modes. In particular, we present the algorithm using a certain set of compactly supported piecewise polynomial test functions which not only have fast decaying spectrum but are designed to yield high accuracy quadrature via the trapezoidal rule, extending the theory from \cite{messenger2020weak}. Pseudocode for the WSINDy\_PDE algorithm is given in \ref{alg:wsindypde}, with sparsity enforced using sequentially-thresholded least squares as in the standard SINDy algorithm \cite{brunton2016discovering} along with Tikhonoff regularization for poorly-conditioned linear systems. 

%We note that the dimensions of both the integration operators $\Vbf^s$ and the data matrix $\Theta(\Ubf)$ have been left ambiguous. 

\subsection{Convolutional Weak Form and Discretization}\label{sec:convweakdisc}

We now restrict to the case of each test function $\psi_k$ being a translation of a reference test function $\psi$, i.e. $\psi_k(x,t) = \psi(\xbf_k-x,t_k-t)$ for some collection of points $\{(\xbf_k,t_k)\}_{k\in [K]} \subset (\Xbf,\tbf)$ (referred to as the {\it query points}). The weak form of the dynamics \eqref{weakform} over the test function basis $(\psi_k)_{k\in [K]}$ then takes the form of a convolution:
\begin{equation}\label{conv_form}
\left(D^{\alpha^0}\psi\right) * u (\xbf_k,t_k) = \sum_{s=1}^S\sum_{j=1}^J \wstar_{(s-1)J+j} \left(D^{\alpha^s}\psi\right) * f_j(u)(\xbf_k,t_k).
\end{equation}
The sign factor $(-1)^{|\alpha^s|}$ appearing in \eqref{weakform} after integrating by parts is eliminated in \eqref{conv_form} due to the sign convention in the integrand of the space-time convolution, which is defined by 
\[\psi*u(x,t) := \int_0^T\int_{\Omega} \psi(x-y,t-s)u(y,s)\,dyds = \lan \psi(x-\cdot, t-\cdot),\ u(\cdot,\cdot)\ran.\]
Construction of the linear system $\bbf = \Gbf\wstar$ as a discretization of the {\it convolutional} weak form \eqref{conv_form} over the query points $\{(\xbf_k,t_k)\}_{k\in[K]}$ can then be carried out efficiently using the FFT as we describe below. 

To relate the continuous and discrete convolutions, we assume that the support of $\psi$ is contained within some rectangular domain
\[\Omega_R := [-b_1,b_1]\times\cdots\times[-b_D,b_D]\times[-b_{D+1},b_{D+1}]\subset \Rbb^D\times\Rbb\]
where $b_d = m_d\Delta x$ for $d\in[D]$ and $b_{D+1} = m_{D+1}\Delta t$. We then define a reference computational grid $(\Ybf,\mathfrak{t})\subset \Rbb^D\times \Rbb$ for $\psi$ centered at the origin and having the same sampling rates $(\Delta x,\Delta t)$ as the data $\Ubf$, where $\Ybf = \Ybf_1\otimes \cdots \otimes \Ybf_D$
for $\Ybf_d = (n\Delta x)_{-m_d\leq n\leq m_d}$ and $\mathfrak{t} = (n\Delta t)_{-m_{D+1}\leq n\leq m_{D+1}}$. In this way $\Ybf$ contains $2m_d+1$ points along each dimension $d\in[D]$, with equal spacing $\Delta x$, and $\mathfrak{t}$ contains $2m_{D+1}+1$ points with equal spacing $\Delta t$. As with $(\Xbf,\tbf)$, points in $(\ybf_k,\mathfrak{t}_k)\in(\Ybf,\mathfrak{t})$ take the form
\[(\ybf_k,\mathfrak{t}_k) = (\Ybf_{k_1,\dots,k_D},\mathfrak{t}_{k_{D+1}})\]
where each index $k_d$ for $d\in[D+1]$ takes values in the range $\{-m_d,\dots,0,\dots,m_d\}$, and for valid indices $k-j$, the two grids $(\Xbf,\tbf)$ and $(\Ybf,\mathfrak{t})$ are related by 
\begin{equation}\label{shiftid}
(\xbf_k-\xbf_j,t_k-t_j) = (\ybf_{k-j},\mathfrak{t}_{k-j}).
\end{equation}
We stress that $(\Ybf,\mathfrak{t})$ is completely defined by the integers $\mbf = (m_d)_{d\in[D+1]}$, specified by the user, and that the values of $\mbf$ have a significant impact on the algorithm. For this reason we develop an automatic selection algorithm for $\mbf$ using spectral properties of the data $\Ubf$ (see Appendix \ref{app:corner}). 

The linear system \eqref{weak_disc} can now be rewritten
\begin{equation}\label{conv_disc}
\begin{dcases} \hspace{1.41cm}\bbf_k = \Psi^0 * \Ubf (\xbf_k,t_k),\\
\Gbf_{k,(s-1)J+j} =  \Psi^s * f_j(\Ubf) (\xbf_k,t_k),\end{dcases}
\end{equation}
where $\Psi^s := D^{\alpha^s}\psi(\Ybf,\mathfrak{t})\Delta x^D\Delta t$ and the factor $\Delta x^D\Delta t$ characterizes the trapezoidal rule. We define the discrete $(D+1)$-dimensional convolution between $\Psi^s$ and $f_j(\Ubf)$ at a point $(\xbf_k,t_k) = (\Xbf_{k_1,\dots,k_D},t_{k_{D+1}})\in (\Xbf,\tbf)$  by
%We define the discrete convolution between two $(D+1)$-dimensional arrays $\Psi$ and $\Ubf$ at the query point $(x_k,t_k) = (\xbf_{k_1,\dots,k_D},\tbf_{k_{D+1}})\in (\xbf,\tbf)$ by 
%\[\Psi*\Ubf(x_k,t_k) := \sum_{j_1=1}^{N_1}\cdots\sum_{j_{D+1}=1}^{N_{D+1}-1} \Psi_{k_1-j_1,\dots,k_{D+1}-j_{D+1}}\Ubf_{j_1,\dots,j_{D+1}},\]
\[\Psi^s*f_j\left(\Ubf\right)(\xbf_k,t_k) := \sum_{\ell_1=1}^{N_1}\cdots\sum_{\ell_{D+1}=1}^{N_{D+1}} \Psi^s_{k_1-\ell_1,\dots,k_{D+1}-\ell_{D+1}} f_j\left(\Ubf_{\ell_1,\dots,\ell_{D+1}}\right),\]
which, substituting the definition of $\Psi^s$,
\begin{align}
&:= \sum_{\ell_1=1}^{N_1}\cdots\sum_{\ell_{D+1}=1}^{N_{D+1}} D^{\alpha^s}\psi\left(\Ybf_{k_1-\ell_1,\dots,k_D-\ell_D},\ \mathfrak{t}_{k_{D+1}-\ell_{D+1}}\right) f_j\left(\Ubf_{\ell_1,\dots,\ell_{D+1}} \right) \Delta x^D\Delta t\\
\intertext{truncating indices appropriately and using \eqref{shiftid},}
\label{convline}&= \sum_{\ell_1=k_1-m_1}^{k_1+m_1}\cdots\sum_{\ell_{D+1}=k_{D+1}-m_{D+1}}^{k_{D+1}+m_{D+1}} D^{\alpha^s}\psi\left(\Ybf_{k_1-\ell_1,\dots,k_D-\ell_D},\ \mathfrak{t}_{k_{D+1}-\ell_{D+1}}\right) f_j\left(\Ubf_{\ell_1,\dots,\ell_{D+1}}\right) \Delta x^D\Delta t\\[10pt]
&= \sum_{\ell_1=k_1-m_1}^{k_1+m_1}\cdots\sum_{\ell_{D+1}=k_{D+1}-m_{D+1}}^{k_{D+1}+m_{D+1}} D^{\alpha^s}\psi\left(\xbf_k-\xbf_\ell,t_k-t_\ell\right) f_j\left(\Ubf_{\ell_1,\dots,\ell_{D+1}}\right) \Delta x^D\Delta t\\[10pt]
&\approx \int_0^T\int_\Omega D^{\alpha^s}\psi(\xbf_k-x,t_k-t) f_j\left(u(x,t)\right) \,dx\,dt.
\end{align}
\subsection{FFT-based Implementation and Complexity for Separable $\psi$}

Convolutions in the linear system \eqref{conv_disc} may be computed rapidly if the reference test function $\psi$ is separable over the given coordinates, i.e.
\[\psi(x,t) =  \phi_1(x_1) \cdots \phi_2(x_D)\phi_{D+1}(t)\]
for univariate functions $(\phi_d)_{d\in[D+1]}$. In this case,
\[D^{\alpha^s}\psi(\Ybf,\mathfrak{t}) = \phi_1^{(\alpha^s_1)}(\Ybf_1)\otimes \dots \otimes\phi_D^{(\alpha^s_D)}(\Ybf_D)\otimes\phi_{D+1}^{(\alpha^s_{D+1})}(\mathfrak{t}),\]
so that only the vectors
\[\phi_d^{(\alpha^s_d)}(\Ybf_d) \in\Rbb^{2m_d+1}, \quad d\in[D] \txt{0.4}{and} \phi_{D+1}^{(\alpha^s_{D+1})}(\mathfrak{t}) \in\Rbb^{2m_{D+1}+1},\]
need to be computed for each $0\leq s \leq S$ and the multi-dimensional arrays $(\Psi^s)_{s=0,\dots,S}$ are never directly constructed. Convolutions can be carried out sequentially in each coordinate\footnote{The technique of exploiting separability in high-dimensional integration is not new (see \cite{pereyra1973efficient} for an early introduction) and is frequently utilized in scientific computing (see \cite{beylkin2005algorithms,harrison2016madness} for examples in computational chemistry). 
}, so that the overall cost of computing each column $\Psi^s*f_j\left(\Ubf\right)$ of $\Gbf$ is
\begin{equation}\label{cfft}
T_{\RNum{1}}(N,n,D) := CN\log(N)\sum_{d=1}^{D+1}N^{D+1-d}\left(N-n+1\right)^{d-1},
\end{equation}
if the computational grid $(\Xbf,\tbf)$ and reference grid $(\Ybf, \mathfrak{t})$ have $N$ and $n\leq N$ points along each of the $D+1$ dimensions, respectively. Here $CN\log(N)$ is the cost of computing the 1D convolution between column vectors $\xbf\in \Rbb^n$ and $\ybf\in \Rbb^N$ using the FFT,
\begin{equation}\label{convfft}
\xbf*\ybf = \CalP\CalF^{-1}\left(\CalF(\xbf^0)\odot\CalF(\ybf)\right),
\end{equation}
where $\xbf^0 = [\,0\cdots 0\ \xbf^T\,]^T \in \Rbb^N$, $\odot$ denotes element-wise multiplication and $\CalP$ projects onto the first $N-n+1$ components. The discrete Fourier transform $\CalF$ is defined 
\begin{align*}
\CalF_k(\ybf) &= \sum_{j=1}^N \ybf_j e^{-2\pi i (j-1)(k-1)}\\
\intertext{with inverse}
\CalF_k^{-1}(\zbf) &= \frac{1}{N}\sum_{j=1}^N \zbf_j e^{2\pi i (j-1)(k-1)}.
\end{align*}
The projection $\CalP$ ensures that the convolution only includes points that correspond to integrating over test functions $\psi$ that are compactly supported in $(\Xbf,\tbf)$, which is necessary for integration by parts to hold in the weak form. The spectra of the test functions $\phi_d^{(\alpha^s_d)}(\Ybf_d)$ can be precomputed and in principle each convolution $\Psi^s*f_j(\Ubf)$ can be carried out in parallel\footnote{For the examples in Section \ref{sec:examples} the walltimes are reported for serial computation of $(\Gbf,\bbf)$.}, making the total cost of the WSINDy Algorithm (\ref{alg:wsindypde}) in the PDE setting equal to \eqref{cfft} (ignoring the cost of the least-squares solves which are negligible in comparison to computing $(\Gbf,\bbf)$). In addition, subsampling reduces the term $(N-n+1)$ in \eqref{cfft} to $(N-n+1)/s$ where $s\geq 1$ is the subsampling rate such that $(N-n+1)/s$ points are kept along each dimension.

\begin{figure}
\begin{tabular}{ccc}
\includegraphics[trim={25 5 42 15},clip,width=0.3\textwidth]{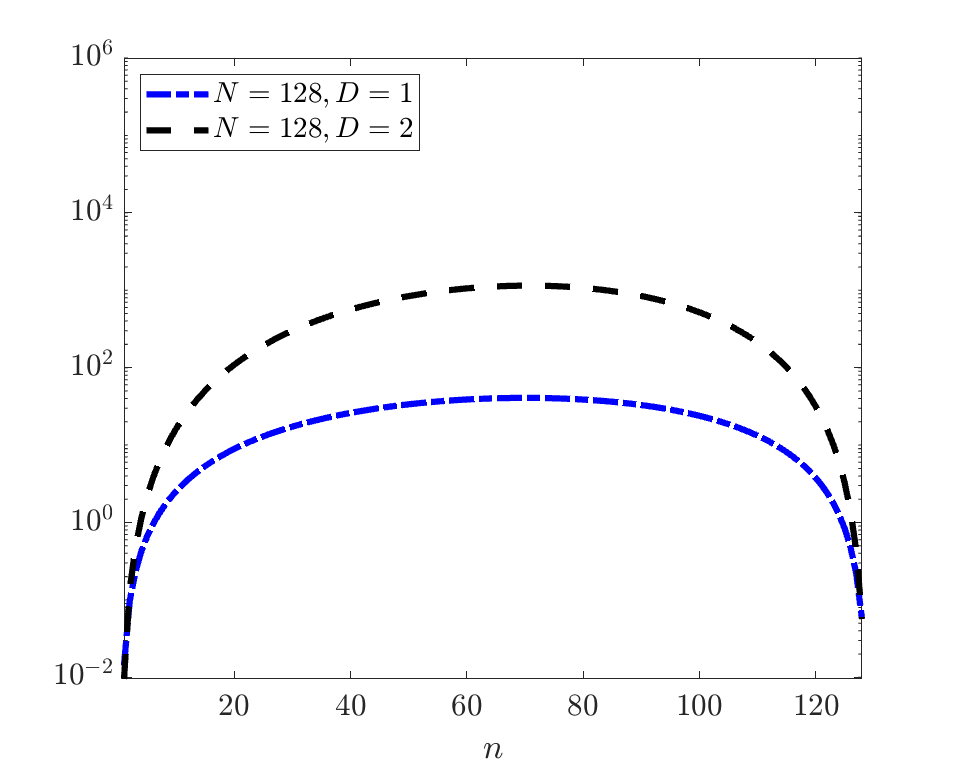} & 
\includegraphics[trim={25 5 41 15},clip,width=0.3\textwidth]{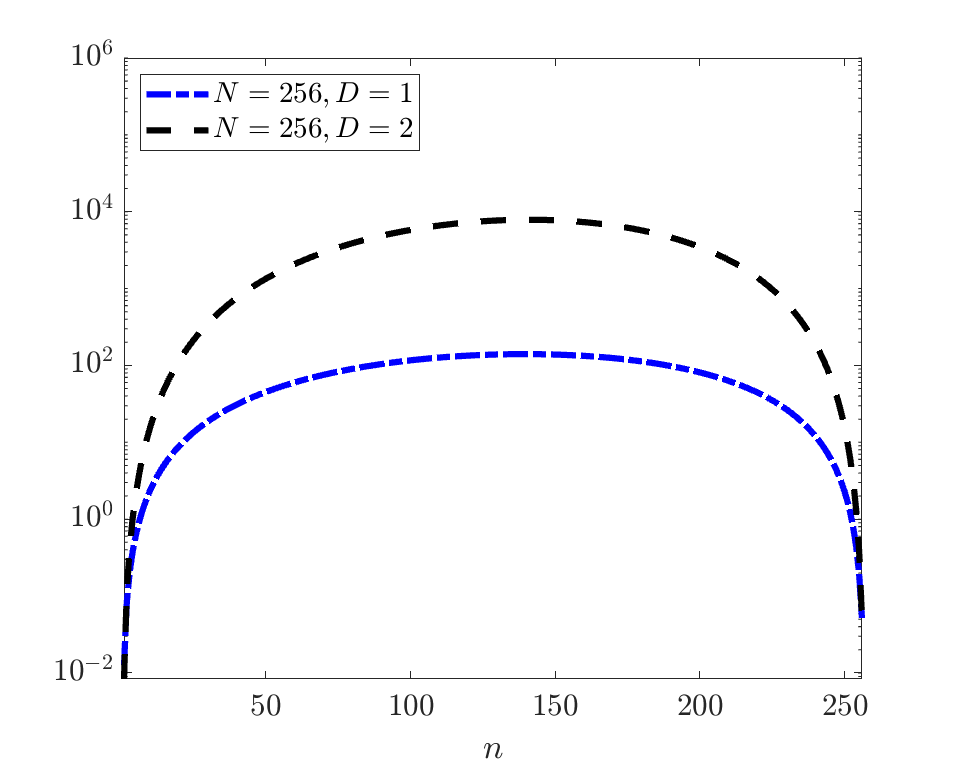} &
\includegraphics[trim={25 5 40 15},clip,width=0.3\textwidth]{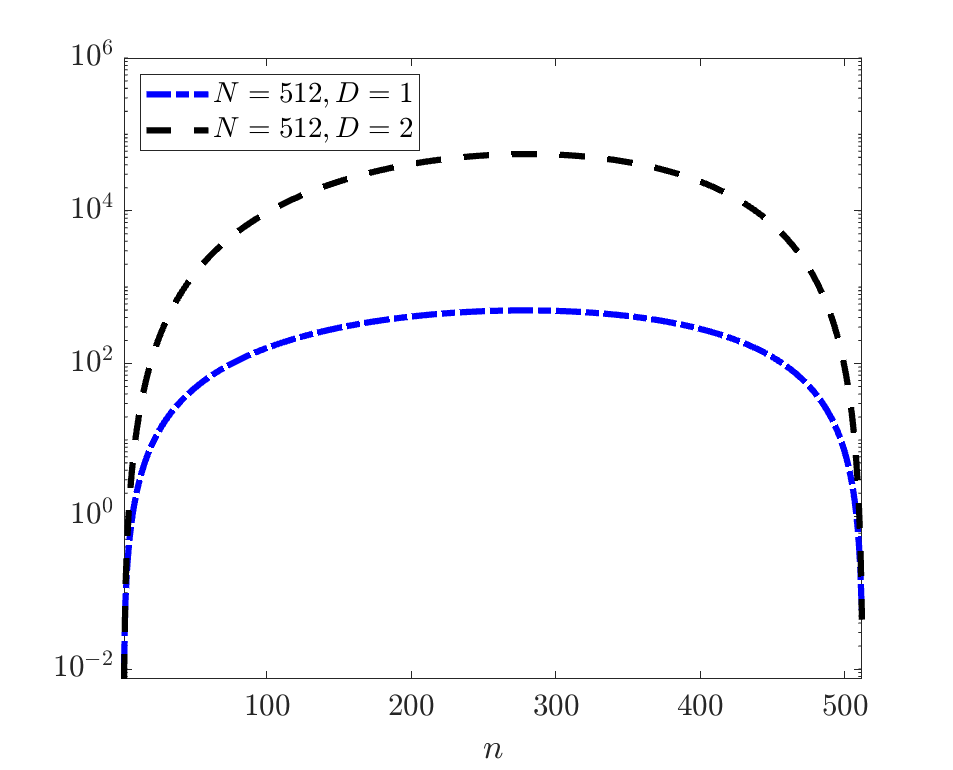} \\ 
\end{tabular}
\caption{Reduction in computational cost of multi-dimensional convolution $\Psi^s*f_j\left(\Ubf\right)$ when $\Psi^s$ and $f_j(\Ubf)$ have $n$ and $N$ points in each of $D+1$ dimensions, respectively. Each plot shows the ratio $T_{\RNum{2}}/T_{\RNum{1}}$ (equations \eqref{cnaive} and \eqref{cfft}), i.e.\ the factor by which the separable FFT-based convolution reduces the cost of the naive convolution, for $D+1 = 2$ and $D+1=3$ space-time dimensions and $n\in[N]$. The right-most plot shows that when $N=512$ and $D+1=3$, the separable FFT-based convolution is $10^4$ times faster for $100\leq n\leq 450$.} 
\label{conv_comp} 
\end{figure}

For most practical combinations of $n$ and $N$, (say $n > N/10$ and $N > 150$) using the FFT and separability provides a considerable reduction in computational cost. See Figure \ref{conv_comp} for a comparison between $T_{\RNum{1}}$ and the naive cost $T_{\RNum{2}}$ of an $(N+1)$-dimensional convolution:
\begin{equation}\label{cnaive}
T_{\RNum{2}}(N,n,D) := (2n^{D+1}-1)(N-n+1)^{D+1}.
\end{equation}
For example, with $n=N/4$ (a typical value) we have $T_{\RNum{2}} = \CalO(N^{2D+2})$ and $T_{\RNum{1}}  =\CalO(N^{D+1}\log(N))$, hence exploiting separability reduces the complexity by a factor of $N^{D+1}/\log(N)$. 

%In addition, the FFT-based convolution reveals a mechanism for ensuring the algorithm's robustness to noise: spectral properties of the test functions can be optimized to damp high-frequency noise. We discuss this further in Section \ref{sec:piecewisepoly}.

%we use separable test functions throughout, although one could easily choose test functions outside of this class. 

\section{WSINDy Algorithm for PDEs and Hyperparameter Selection}\label{sec:alg}

WSINDy for PDE discovery is given in Algorithm \ref{alg:wsindypde}, where the user must specify each of the hyperparameters in Table \ref{hypparam}. The key pieces of the algorithm are (i) the choice of reference test function $\psi$, (ii) the method of a sparsification, (iii) the method of regularization, (iv) selection of convolution query points $\{(\xbf_k,t_k)\}_{k\in K}$, and (v) the model library. At first glance, the number of hyperparameters is quite large. We now discuss several simplifications that either reduce the number of hyperparameters or provide methods of choosing them automatically. In Section \ref{sec:piecewisepoly} we discuss connections between the convolutional weak form and spectral properties of $\psi$ that determine the scheme's robustness to noise and inform the selection of test function hyperparameters. In Section \ref{sec:sparse} we introduce a {\it modified} sequential-thresholding least-squares algorithm (MSTLS) which includes automatic selection of the threshold $\lambda$ and allows for PDE discovery from large libraries. In Section \ref{sec:rescale} we describe how scale-invariance of the PDE is used to rescale the data and coordinates in order to regularize the model recovery problem in the case of poorly-scaled data. In Sections \ref{sec:query} and \ref{sec:library} we briefly discuss selection of query points and an appropriate model library, however these components of the algorithm will be investigated more thoroughly in future research.

\subsection{Selecting a Reference Test Function $\psi$}\label{sec:piecewisepoly}

\subsubsection{Convolutional Weak Form and Fourier Analysis} Computation of $\Gbf$ and $\bbf$ in \eqref{conv_disc} with $\psi$ separable requires the selection of appropriate 1D coordinate test functions $(\phi_d)_{d\in[D+1]}$. Computing convolutions using the FFT \eqref{convfft} suggests a mechanism for choosing appropriate test functions. Define the Fourier coefficients of a function $u\in L^2([0,T])$ by 
\[\widehat{u}(k) = \frac{1}{\sqrt{T}}\int_0^T u(t) e^{-\frac{2\pi ik}{T}t}\,dt, \qquad k\in \Zbb.\]
Consider data $\Ubf = u(\tbf) +\ep \in \Rbb^N$ for a $T$-periodic function $u$, $\tbf_k = k\frac{T}{N} = k\Delta t$, and white noise $\ep\sim \CalN(0,\sigma^2 \Ibf)$. The discrete Fourier transform of the noise $\CalF(\ep) := \ep_R+i\ep_I$ is then distributed $\ep_R,\ep_I\sim \CalN(0,(N\sigma^2/2) \Ibf)$. In addition, there exist constants $C>0$ and $\ell>1/2$ such that $|\widehat{u}_k|\leq C|k|^{-\ell}$ for each $k\in \Zbb$. There then exists a noise-dominated region of the spectrum $\CalF(\Ubf)$ determined by the noise-to-signal ratio
\[NSR_k:= \Ebb\left[\frac{|\CalF_k(\ep)|^2}{|\CalF_k(u(\xbf))|^2}\right] = \frac{N\sigma^2}{|\CalF_k(u(\xbf))|^2} \approx \frac{T\sigma^2}{N|\widehat{u}(k)|^2}  \geq \frac{1}{C^2}\Delta t\sigma^2 k^{2\ell},\]
where `$\approx$' corresponds to omitting the aliasing error. For $NSR_k\geq 1$ the $k$th Fourier mode is by definition noise-dominated, which corresponds to wavenumbers
\begin{equation}\label{cornerptdef}
|k| \geq k^*\approx  \left(\frac{C}{\sigma \sqrt{\Delta t}}\right)^{1/\ell}.
\end{equation}
If the critical wavenumber $k^*$ between the noise dominated ($NSR_k\geq 1$) and signal-dominated ($NSR_k\leq 1$) modes can be estimated from the dataset $\Ubf$,
then it is possible to design test functions $\psi$ such that the noise-dominated region of $\CalF(\Ubf)$ lies in the tail of $\widehat{\psi}$. The convolutional weak form (\ref{conv_disc}) can then be interpreted as an approximate low-pass filter on the noisy dataset, offering robustness to noise without altering the frequency content of the data\footnote{This is in contrast to explicit data-denoising, where a filter is applied to the dataset prior to system identification and may fundamentally alter the underlying clean data. The implicit filtering of the convolutional weak form is made explicit by the FFT-based implementation \eqref{convfft}.}.

In summary, spectral decay properties of the reference test function $\psi$ serve to damp high-frequency noise in the convolutional weak form, which acts together with the natural variance-reducing effect of integration, as described in \cite{gurevich2019robust}, to allow for quantification and control of the scheme's robustness to noise. Specifically, coordinate test functions $\phi_d$ with wide support in real space (larger $m_d$) will reduce more variance, but will have a faster-decaying spectrum $\widehat{\phi}_d$, so that \textit{signal}-dominated modes may not be resolved, leading to model misidentification. On the other hand, if $\phi_d$ decays too swiftly in real space (smaller $m_d$), then the spectrum $\widehat{\phi}_d$ will decay more slowly and may put too much weight on \textit{noise}-dominated frequencies. In addition, smaller $m_d$ may not sufficiently reduce variance. A balance must be struck between (a) effectively reducing variance, which is ultimately determined by the decay of $\psi$ in physical space, and (b) resolving the underlying dynamics, determined by the decay of $\widehat{\psi}$ in Fourier space.

\subsubsection{Piecewise-Polynomial Test Functions} 
\begin{figure}
\includegraphics[trim={100 30 100 15},clip,width=0.95\textwidth]{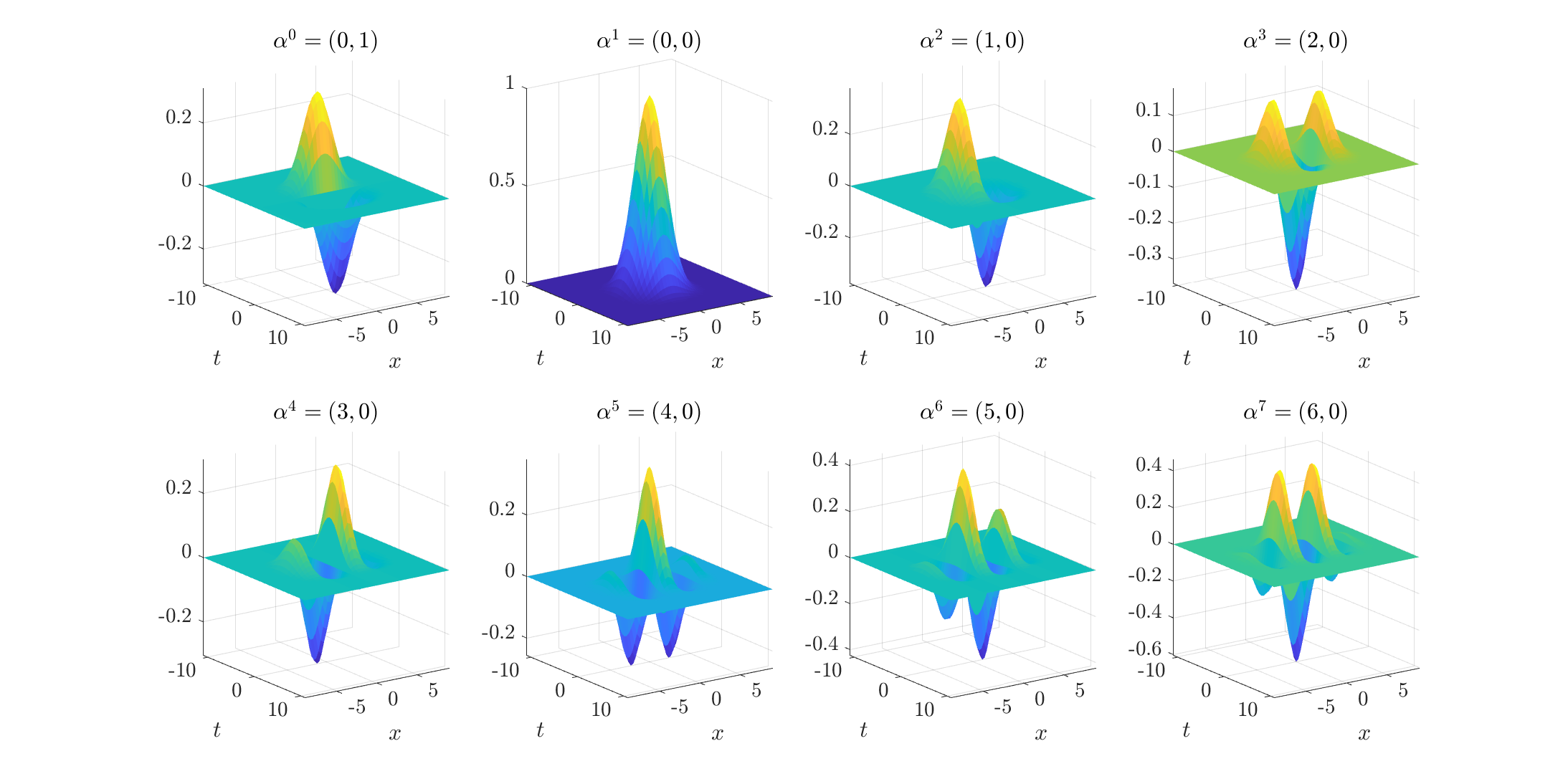} 
\caption{Plots of reference test function $\psi$ and partial derivatives $D^{\alpha^s}\psi$ used for identification of the Kuramoto-Sivashinsky equation. The upper left plot shows $\partial_t\psi$, the bottom right shows $\partial_x^6 \psi$. See Tables \ref{pdetable}-\ref{spectable} for more details.} 
\label{basis_fig} 
\end{figure}

Many test functions achieve the necessary balance between decay in real space and decay in Fourier space in order to offer both variance reduction and resolution of signal-dominated modes (defined by \eqref{cornerptdef}). For simplicity, in this article we use the same test function space used in the ODE setting \cite{messenger2020weak} and leave an investigation of the performance of different test functions to future work. Define $\CalS$ to be the space of functions 
\begin{equation}
\phi(v)=\begin{cases}
C(v-a)^{p}(b-v)^{q} & a<v<b,\\
0 & \text{otherwise},
\end{cases}\label{testfcn}
\end{equation}
where $p,q\geq 1$ and $v$ is a variable in time or space. The normalization 
\[
C=\frac{1}{p^{p}q^{q}}\left(\frac{p+q}{b-a}\right)^{p+q}
\]
ensures that $\left\Vert \phi\right\Vert _{\infty}=1$. Functions
$\phi\in\mathcal{S}$ are non-negative, unimodal, compactly-supported in $[a,b]$, and have $\lfloor\min\{p,q\}\rfloor$ weak derivatives. Larger $p$ and $q$ imply faster decay towards the endpoints $(a,b)$ and for $p=q$ we refer to $p$ as the {\it degree} of $\phi$. See Figure \ref{basis_fig} for a visualization of $\psi$ and partial derivatives $D^{\alpha^s}\psi$ constructed from tensor products of functions from $\CalS$. In addition to having nice integration properties combined with the trapezoidal rule (see Lemma 1 of \cite{messenger2020weak}), $(a,b,p,q)$ can be chosen to localize $\widehat{\phi}$ around signal-dominated frequencies in $\CalF(\Ubf)$ using that $|\widehat{\phi}(k)|=o\left(|k|^{-\lfloor\min\{p,q\}\rfloor-1/2}\right)$.

To assemble the reference test function $\psi$ from one-dimensional test functions $(\phi_d)_{d\in[D+1]}\subset \CalS$ along each coordinate, we must determine the parameters $(a_d,b_d,p_d,q_d)$ in the formula \eqref{testfcn} for each $\phi_d$. Letting $p_d=q_d$, so that each $\phi_d$ is symmetric, and centering $(\Ybf,\mathfrak{t})$ at the origin, we  have that each $\phi_d$ is supported on $[a_d,b_d] = [-b_d,b_d]$ where $b_d = m_d\Delta x$ for $d\in[D]$ and $b_{D+1} = m_{D+1}\Delta t$, so that only $\{(m_d,p_d)\}_{d\in[D+1]}$ need to be specified. The vectors $(\phi_d^{(\alpha^s_d)}(\Ybf_d))_{0\leq s\leq S}$ can be computed from an analogous function $\overline{\phi}_{p_d}$ with support $[-1,1]$, 
\[\overline{\phi}_{p_d}(v) := \begin{dcases} (1-v^2)^{p_d}, & -1<v<1 \\ 0, &\text{otherwise}, \end{dcases}\]
using
\[\phi_d^{(\alpha^s_d)}(\Ybf_d) = \frac{1}{b_d^{\alpha^s_d}}\overline{\phi}_{p_d}^{(\alpha^s_d)}\left(\frac{\Ybf_d}{b_d}\right) = \frac{1}{(m_d\Delta)^{\alpha^s_d}}\overline{\phi}_{p_d}^{(\alpha^s_d)}\left(\nbf_d\right),\]
where the scaled grid $\nbf_d$ is defined $\nbf_d := (n/m_d)_{-m_d\leq n\leq m_d}$ and $\Delta\in\{\Delta x,\Delta t\}$.  

The discrete support lengths $\mbf = (m_d)_{d\in [D+1]}$ and degrees $\pbf = (p_d)_{d\in[D+1]}$ determine the smoothness of $\psi$, as well as its decay in real and in Fourier space, hence are critical to the method's performance. The degrees $\pbf$ can be chosen from $\mbf$ to ensure necessary smoothness and decay in real space using
\begin{equation}\label{realdecay}
p_d = \min\left\{p\geq \overline{\alpha}_d +1 \txt{0.2}{:} \overline{\phi}_{p}\left(1-\frac{1}{m_d}\right) \leq \tau\right\},
\end{equation}
where $\overline{\alpha}_d := \max_{0\leq s\leq S}(\alpha^s_d)$ is the maximum derivative along the $d$th coordinate and $\tau$ is a chosen decay tolerance. In this way $\phi_d$ decays to $\tau$ at the first interior gridpoint of its support, which controls the integration error, and $\phi_d \in C^{\overline{\alpha}_d}(\Rbb)$ so that $\psi$ is smooth enough to integrate by parts as many times as required by the multi-index set $\pmb{\alpha}$. Altogether, the steps for arriving at the test function values on the reference grid $(\phi_d^{(\alpha^s_d)}(\Ybf_d))_{0\leq s\leq S}$ are contained in Algorithm \ref{alg:gettestfcn}. 

In the examples below, we set $\tau = 10^{-10}$ throughout\footnote{WSINDy appears not to be particularly sensitive to $\tau$, similar results were obtained for $\tau = 10^{-6}, 10^{-10}, 10^{-16}$.} and we use the method introduced in Appendix \ref{app:corner} to choose $\mbf$, which involves estimating the critical wavenumber $k^*$ (defined in \eqref{cornerptdef}) between noise-dominated and signal-dominated modes of $\CalF(\Ubf)$. We also simplify things by choosing the same coordinate test function for all spatial coordinates, $\phi_1=\phi_2=\cdots=\phi_D :=\phi_x$ and $\phi_{D+1} = \phi_t$, where $\phi_x$ has degree $p_x$ and support $m_x$ and $\phi_t$ has degree $p_t$ and support $m_t$ (recall $x$ is a sub-index on $\phi_x$ and not a partial derivative). This convention is used in the following sections.

%\begin{algorithm}
%	\caption{$\Cbf =  \textbf{get\_scaling\_fcn\_vals}\left(\Delta x,m,p,\overline{\alpha}\right)$:}
%	\label{alg:gettestfcn}
%	\begin{algorithmic}[1]		
%				\STATE Initialize	$\Cbf \in \Rbb^{(\overline{\alpha}+1)\times 2m+1} = \textbf{0}\in \Rbb^{(\overline{\alpha}+1)\times(2m+1)}$
%				\STATE Initialize $\Lbf = \textbf{0}\in \Rbb^{(\overline{\alpha}+1)\times(m+1)}$,$\Rbf = \textbf{0}\in \Rbb^{(\overline{\alpha}+1)\times(m+1)}$, $\qbf = \textbf{1}\in\Rbb^{\overline{\alpha}+1}$
%				\FOR{$j=1:m+1$}
%							\STATE Set $\Lbf_{1,j} = (1+(j-1)/m)^{p-\overline{\alpha}}$
%							\STATE Set $\Rbf_{1,j} = (-1)^{\overline{\alpha}}(1-(j-1)/m)^{p-\overline{\alpha}}$
%							\FOR{$k=1:\overline{\alpha}$}
%								\STATE Set $\Lbf_{k+1,j} = (1+(j-1)/m)\Lbf_{k,j}$
%								\STATE Set $\Rbf_{k+1,j} = (-1)(1-(j-1)/m)\Rbf_{k,j}$
%							\ENDFOR
%				\ENDFOR
%				\STATE Set $\qbf_{2} = p$
%				\FOR{$\ell = 1:\overline{\alpha}-1$}
%						\STATE $\qbf_{\ell+2} = (p-\ell) \qbf_{\ell+1}$
%				\ENDFOR
%				\FOR{$j=1:m+1$}
%						\FOR{$\ell=1:\overline{\alpha}+1$}
%								\FOR{$q=0:\ell$}
%								 		\STATE Set $\Cbf_{\ell,m+j}  = \Cbf_{\ell,m+j} + {\ell \choose q}\qbf_{\ell-q}\Lbf_{\ell-q,j}\qbf_{q}\Rbf_{q,j}$				
%								\ENDFOR
%						 		\STATE Set $\Cbf_{\ell,m+j}  = \Cbf_{\ell,m+j} + {\ell \choose q}\qbf_{\ell-q}\Lbf_{\ell-q,j}\qbf_{q}\Rbf_{q,j}$				

%												

%\end{algorithmic}
%\end{algorithm}

\begin{algorithm}
	\caption{$(\phi_d^{(\alpha^s_d)}(\Ybf_d))_{0\leq s\leq S}=  \textbf{get\_test\_fcns}\left(m_d, \tau;\ \Xbf_d, \boldsymbol{\alpha}\right)$:}
	\label{alg:gettestfcn}
	\begin{algorithmic}[1]		
	\STATE $N_d$ = length($\Xbf_d$)
	\STATE $\Delta x$ = gridwidth($\Xbf_d$)
						\IF{$m_d > \frac{N_d-1}{2}$ or $m_d\leq 1$}
								\STATE return (``ERROR: invalid support size $m_d$'')
								\STATE BREAK
						\ENDIF
				\STATE Set $\overline{\alpha}_d = \max_{0\leq s\leq S} (\alpha^s_d)$
				\STATE Solve $p_d = \min\left\{p\geq \overline{\alpha}_d +1\txt{0.2}{:} \overline{\phi}_{p}\left(1-\frac{1}{m_d}\right) \leq \tau\right\}$
				\STATE Initialize $\Abf = \textbf{0}\in\Rbb^{(S+1)\times (2m_d+1)}$ 
				\STATE Set $\nbf_d := (n/m_d)_{-m_d\leq n\leq m_d}$
				\FOR{$s=0:S$}
						\STATE Compute analytical order-$(\alpha^s_d)$ derivatives $\Abf_s = \overline{\phi}_{p_d}^{(\alpha^s_d)}(\nbf_d)$ \\
\STATE Set $\phi_d^{(\alpha^s_d)}(\Ybf_d) = \frac{1}{(m_d\Delta x)^{\alpha^s_d}}\Abf_s$
				\ENDFOR
\end{algorithmic}
\end{algorithm}

\subsection{Sparsification}\label{sec:sparse}
To enforce a sparse solution we present a \textit{modified sequential-thresholding least-squares} algorithm MSTLS$(\Gbf,\bbf;\,\lambda)$, defined in \eqref{MSTLS1}, which accounts for terms that are outside of the dominant balance physics of the data, as determined by the left-hand side $\bbf$, as well as terms with small coefficients. We then utilize the loss function
\begin{equation}\label{lossfcn}
\CalL(\lambda) = \frac{\nrm{\Gbf(\wbf^\lambda-\wbf^{LS})}_2}{\nrm{\Gbf\wbf^{LS}}_2}+\frac{\#\CalI^\lambda}{SJ}
\end{equation}
to select an optimal threshold $\widehat{\lambda}$, where $\wbf^\lambda$ is the output of MSTLS$(\Gbf,\bbf;\,\lambda)$, $\#\CalI^\lambda$ is the cardinality of the index set $\CalI^\lambda := \{j\,:\,\wbf^\lambda_j\neq 0\}$ of non-zero coefficients, $\wbf^{LS} := \left(\Gbf^T \Gbf \right)^{-1} \Gbf ^T \bbf $ is the least squares solution, and $SJ$ is the total number of terms in the library ($S$ differential operators and $J$ nonlinear functions). The two terms in $\CalL$ penalize (i) the distance between $\Gbf \wbf^{LS}$ (the projection of $\bbf$ onto the range of $\Gbf$) and $\Gbf \wbf^\lambda$ (the projection of $\bbf$ onto the range of the restriction $\Gbf_{\CalI^\lambda}$) and (ii) the number of nonzero terms in the resulting model, respectively, with normalization ensuring that $\CalL(0) = \CalL(\infty) = 1$. 

The $\text{MSTLS}(\Gbf,\bbf;\,\lambda)$ iteration is as follows. For a given $\lambda\geq 0$, define the set of lower bounds $L^\lambda$ and upper bound $U^\lambda$ by 
\[\begin{dcases} L_j^\lambda =  \lambda\max\left\{1,\ \frac{\nrm{\bbf}}{\nrm{\Gbf_j}}\right\}\\
U_j^\lambda =  \frac{1}{\lambda}\min\left\{1,\ \frac{\nrm{\bbf}}{\nrm{\Gbf_j}}\right\}\end{dcases}, \qquad 1\leq j\leq SJ.\]
Then with $\wbf^0 = \wbf^{LS}$, define the iterates 
\begin{equation}\label{MSTLS1}
\begin{dcases} \hspace{0.43cm}\CalI^\ell = \{1\leq j\leq SJ\ :\ L^\lambda_j\leq|\wbf^\ell_j|\leq U^\lambda_j\} \\
\wbf^{\ell+1} = \argmin_{\supp{\wbf}\subset \CalI^\ell} \nrm{ \Gbf  \wbf-\bbf}_2^2.\end{dcases}
\end{equation}
The constraint $L^\lambda_j\leq|\wbf^\ell_j|\leq U^\lambda_j$ is clearly more restrictive than standard sequential thresholding, but it enforces two desired qualities of the model: (i) that the coefficients $\wbf^\lambda$ do not differ too much from 1, since 1 is the coefficient of the ``evolution'' term $D^{\alpha^0}u$ (assumed known), and (ii) that the ratio $\nrm{\wbf_j \Gbf_j}_2/\nrm{ \bbf }_2$ lies in $[\lambda,\lambda^{-1}]$ enforcing an empirical dominant balance rule (e.g.\ $\lambda=0.01$ allows terms in the model to be at most two orders of magnitude from $D^{\alpha^s}u$). The overall sparsification algorithm $\text{MSTLS}(\Gbf,\bbf;\, \CalL,\pmb{\lambda})$ is
\begin{equation}\label{MSTLS2}
\begin{dcases} 
\hspace{0.1cm} \widehat{\lambda} = \min\left\{\lambda\in \pmb{\lambda} \ :\ \CalL(\lambda) = \min_{\lambda\in \pmb{\lambda}} \CalL(\lambda)\right\}\\
\widehat{\wbf}  =\text{MSTLS}(\Gbf,\bbf;\,\widehat{\lambda}),\end{dcases}
\end{equation}
where $\pmb{\lambda}$ is a finite set of candidate thresholds\footnote{Other methods of minimizing $\CalL$ can be used, however minimizers are not unique (there exists a set of minimizers - see Figure \ref{loss}). Our approach is efficient and returns the minimizer $\widehat{\lambda}$ which has the useful characterization of defining the thresholds $\lambda$ that result in overfitting.}. The learned threshold $\widehat{\lambda}$ is the smallest minimizer of $\CalL$ over the range $\pmb{\lambda}$ and hence marks the boundary between identification and misidentification of the minimum-cost model, such that $\{\lambda \in \pmb{\lambda}\,:\, \lambda <\widehat{\lambda}\}$ results in overfitting. A similar learning method for $\widehat{\lambda}$ combining STLS and Tikhonov regularization (or {\it ridge regression}) was developed in \cite{rudy2017data}. We have found that our approach of combining $\text{MSTLS}(\Gbf,\bbf;\, \CalL,\pmb{\lambda})$ with rescaling, as introduced in the next section, regularizes the sparse regression problem in the case of large model libraries without adding hyperparameters\footnote{Tikhonov regularization involves solving $\what = \argmin_\wbf\nrm{\Gbf \wbf -\bbf}_2^2+\gamma^2\nrm{\wbf}_2^2$} and definitely deserves further study. 

\subsection{Regularization through Scale Invariance}\label{sec:rescale}

Construction of the linear system $\bbf = \Gbf \wbf$ involves taking (nonlinear) transformations of the data $f_j(\Ubf)$ and then integrating against $D^{\alpha^s}\psi$, which oscillates for large $|\alpha^s|$. This can lead to a large condition number $\kappa(\Gbf)$ and prevent accurate inference of the true model coefficients $\wstar$, especially when the underlying data is poorly scaled\footnote{A common remedy for this is to scale $\Gbf$ to have columns of unit 2-norm, however this has no connection with the underlying physics.}. Often characteristic scales effect the dynamics in nontrivial ways such that naively rescaling the data leads to inference of incorrect model coefficients. For the inviscid Burgers and KdV data below, the amplitude of the data determines the wavespeed, but with $\Ubf = \CalO(10^3)$ identification of the term $\partial_x(u^2)$ from a large library of polynomial nonlinearities is ill-conditioned. To overcome this we propose to rescale the underlying coordinates to achieve low condition number using scale invariance of the PDE. 

If $u$ solves \eqref{diffform}, then for any $\gamma_x,\gamma_t,\gamma_u>0$, the function
\[\tilde{u}(\tilde{x},\tilde{t}) := \gamma_u\, u\left(\frac{\tilde{x}}{\gamma_x},\frac{\tilde{t}}{\gamma_t}\right) := \gamma_u\, u(x,t)\]
solves, 
\[\tilde{D}^{\alpha^0}\tilde{u} = \sum_{s=1}^S\sum_{j=1}^J \tilde{\wbf}_{(s-1)J+j}\tilde{D}^{\alpha^s}\tilde{f}_j(\tilde{u})\]
where $\tilde{D}^{\alpha^s}$ denotes differentiation with respect to $(\tilde{x},\tilde{t}) = (\gamma_xx,\gamma_tt)$
and $\tilde{f}_j(\tilde{u})= f_j(\tilde{u}) = \gamma_u^{\beta_j} f_j(u)$ for homogeneous functions with power $\beta_j$ and $\tilde{f}_j(\tilde{u})= f_j\left(\frac{\tilde{u}}{\gamma_u}\right) = f_j(u)$ otherwise (in the latter case we set $\beta_j=0$). The linear system in the rescaled coordinates $\tilde{\bbf} = \tilde{\Gbf}\tilde{\wbf}$ is constructed by discretizing the convolutional weak form as before but with a reference test function $\tilde{\psi}$ on the rescaled grid $\tilde{\Omega}_R$. We recover the coefficients $\what$ at the original scales by setting $\what = \Mbf \tilde{\wbf}$, where $\Mbf = \text{diag}\left({\pmb{\mu}}\right)$ is the diagonal matrix with entries
\[\mu_{(s-1)J+j} :=  \gamma_u^{-(\beta_j-1)}\,\gamma_x^{\sum_{d=1}^D(\alpha^s_d-\alpha^0_d)}\ \gamma_t^{(\alpha^s_{D+1}-\alpha_{D+1}^0)}.\]
To choose the scales $\gamma_x,\gamma_t$ for $\tilde{\phi}\in \CalS$ we note that in the coordinates $(\tilde{x},\tilde{t})$, for even derivatives $\alpha_x$ and $\alpha_t$ we have
\[\nrm{\tilde{\phi}_x^{(\alpha_x)}}_\infty = \frac{{p_x \choose \frac{\alpha_x}{2}}\alpha_x !}{m_x^{\alpha_x}(\gamma_x\Delta x)^{\alpha_x}}, \quad \nrm{\tilde{\phi}_t^{(\alpha_t)}}_\infty = \frac{{p_t \choose \frac{\alpha_t}{2}}\alpha_t !}{m_t^{\alpha_t}(\gamma_t\Delta t)^{\alpha_t}}\]
and so setting
\begin{equation}\label{scalext}
\gamma_x = \frac{1}{m_x \Delta x}\left({p_x \choose \frac{\overline{\alpha}_x}{2}}\overline{\alpha}_x !\right)^{1/\overline{\alpha}_x},\quad \gamma_t = \frac{1}{m_t \Delta t}\left({p_t \choose \frac{\overline{\alpha}_t}{2}}\overline{\alpha}_t !\right)^{1/\overline{\alpha}_t},
\end{equation}
where $\overline{\alpha}_x$ and $\overline{\alpha}_t$ are the maximum spatial and temporal derivative appearing in the library, ensures that\footnote{Here $\nrm{\Psi^s}_{1'}$ is the 1-norm of $\Psi^s$ streched into a column vector (i.e.\ the trapezoidal-rule approximation of $\int_{\Omega_R} |D^{\alpha^s}\psi|\,dxdt$).} 
\[\max_s\nrm{\Psi^s}_{1'}\leq \max_s\nrm{D^{\alpha^s}\tilde{\psi}}_\infty |\tilde{\Omega}_R|\leq |\tilde{\Omega}_R|.\]
We then set $\gamma_u$ according to the fastest growing term $f_j$ in the library to get an approximate uniform bound on the columns of the scaled Gram matrix $\tilde{\Gbf}$. Since the examples below all use monomials, we let
\begin{equation}\label{scaleu}
\gamma_u = \nrm{\frac{\Ubf^{\overline{\beta}}}{\nrm{\Ubf}_F}}_{F}^{-1/\overline{\beta}}
\end{equation}
where $\overline{\beta} = \max_j\beta_j$, so that $\nrm{\tilde{f}_{j}(\Ubf)}_F \approx \nrm{\Ubf}_F$, and $\nrm{\cdot}_F$ is the Frobenius norm. From Young's inequality for convolutions, we then have 
\[\nrm{\tilde{\Gbf}_{(s-1)J+j}}_F= \nrm{\Psi^s*\tilde{f}_j(\Ubf)}_F\leq \nrm{\Psi^s}_{1'}\nrm{\tilde{f}_j(\Ubf)}_F \approx |\tilde{\Omega}_R| \nrm{\Ubf}_F.\]
Similar scales $\gamma_x,\gamma_t,\gamma_u$ can be chosen for different model libraries and reference test functions. In the examples below we rescale the data and coordinates according to \eqref{scalext} and \eqref{scaleu}, which results in a low condition number $\kappa(\tilde{\Gbf})$ (see Table \ref{spectable}). Throughout what follows, quantities defined over scaled coordinates will be denoted by tildes.

\subsection{Query Points and Subsampling}\label{sec:query}

Placement of $\{(\xbf_k,t_k)\}_{k\in[K]}$ determines which regions of the observed data will most influence the recovered model\footnote{Note that the projection operation in \eqref{convfft} restricts the admissable set of query points to those for which $\psi(\xbf_k-x,t_k-t)$ is compactly supported within $\Omega\times[0,T]$, which is necessary for integration by parts to be valid.}. In WSINDy for ODEs (\cite{messenger2020weak}), an adaptive algorithm was designed for placement of test functions near steep gradients along the trajectory. Improvements in this direction in the PDE setting are a topic of active research, however, for simplicity in this article we uniformly subsample $\{(\xbf_k,t_k)\}_{k\in[K]}$ from $(\Xbf,\tbf)$ using subsampling frequencies $\sbf = (s_1,\dots, s_{D+1})$ along each coordinate, specified by the user. That is, along each one-dimensional grid $\Xbf_d$, $\lfloor\frac{N_d-2m_d}{s_d}\rfloor$ points are selected with uniform spacing $s_d\Delta x$ for $d\in[D]$ and $s_{D+1}\Delta t$ for $d=D+1$. This results in a $(D+1)$-dimensional coarse grid with dimensions $\lfloor\frac{N_1-2m_1}{s_1}\rfloor\times\cdots\times \lfloor\frac{N_{D+1}-2m_{D+1}}{s_{D+1}}\rfloor$, which determines the number of query points
\begin{equation}\label{numK}
K = \prod_{d=1}^{D+1} \left\lfloor\frac{N_d-2m_d}{s_d}\right\rfloor.
\end{equation}

\subsection{Model Library}\label{sec:library}

 The model library is determined by the nonlinear functions $(f_j)_{j\in[J]}$ and the partial derivative indices $\pmb{\alpha}$ and is crucial to the well-posedness of the recovery problem. In the examples below we choose $(f_j)_{j\in[J]}$ to be polynomials or trigonometric functions as these sets are dense in many relevant function spaces. For simplicity in this work we choose $\pmb{\alpha}$ without cross-terms (e.g.\ $\frac{\partial^2}{\partial {x_1} \partial {x_2}}$ is omitted), however including these terms in the library does not have a significant impact on the results below.

\begin{algorithm}
	\caption{$(\what, \, \hat{\lambda}) =$ \textbf{WSINDy}$((f_j)_{j\in [J]},\,\boldsymbol{\alpha},\mbf,\sbf,\,\pmb{\lambda},\tau;\ \Ubf, (\Xbf,\tbf))$:}
	\label{alg:wsindypde}
	\begin{algorithmic}[1]
		\FOR{$d=1:D+1$}
				\STATE Compute $(\phi_d^{(\alpha^s_d)}(\Ybf_d))_{0\leq s\leq S}=  \textbf{get\_test\_fcns}\left(m_d, \tau;\ \Xbf_d, \boldsymbol{\alpha}\right)$ using Algorithm \ref{alg:gettestfcn}
%				\STATE Compute one-dimensional test functions $(\phi_d^{(\alpha_d^s)}(\xbf_d))_{d\in [D+1]}$ from support sizes $\abf = (a_1, \dots, a_{D+1})$ such that $\Psi^s = \phi_d^{(\alpha_d^s)}(\xbf_d) \otimes \cdots \otimes \phi_{D+1}^{(\alpha_{D+1}^s)}(\tbf) \, (\Delta x)^D\Delta t$;
		\ENDFOR 
		\STATE Compute scales $\{\gamma_u,(\gamma_d)_{d=1}^{D+1}\}$ and scale matrix $\Mbf = \text{diag}(\pmb{\mu})$ using \eqref{scalext} and \eqref{scaleu}
		\STATE Subsample query points $\{(\xbf_k,t_k)\}_{k\in[K]}\subset (\Xbf,\tbf)$ using subsampling frequencies $\sbf = (s_1,s_2,\dots, s_{D+1})$;
		\STATE Compute left-hand side $\tilde{\bbf} = \tilde{\Psi}^0 * \tilde{\Ubf}$ over $\{(\xbf_k,t_k)\}_{k\in[K]}$ using FFT and separability  of $\psi$;
		\FOR{$j=1:J$}
				\STATE Compute $\tilde{f}_j(\tilde{\Ubf})$;
				\FOR{$s=1:S$}
						\STATE Compute column $(s-1)J+j$ of Gram matrix $\tilde{\Gbf}_{:,(s-1)J+j} =  \tilde{\Psi}^s * \tilde{f}_j(\tilde{\Ubf})$ over $\{(\xbf_k,t_k)\}_{k\in[K]}$ using FFT and separability of $\psi$
 			 \ENDFOR
		\ENDFOR
		\STATE $(\widehat{\wbf},\widehat{\lambda}) = \text{MSTLS}(\tilde{\Gbf},\tilde{\bbf};\, \CalL, \pmb{\lambda})$
\end{algorithmic}
\end{algorithm}

\begin{table}
\centering
\begin{tabular}{|c|c|c|}
\hline
Hyperparameter & Domain & Description \\\hline
$(f_j)_{j\in[J]}$ & $BV_{loc}(\Rbb)$ & trial function library\\ \hline
$\pmb{\alpha} = (\alpha_s)_{s=0,\dots,S}$ & $\Nbb^{(S+1) \times (D+1)}$ & partial derivative multi-indices \\ \hline
$\mbf = (m_d)_{d\in[D+1]}$ & $\Nbb^{D+1}$ & discrete support lengths of 1D test functions $(\phi_d)_{d\in[D+1]}$ \\ \hline
$\sbf = (s_d)_{d\in[D+1]}$ & $\Nbb^{D+1}$ & subsampling frequencies for query points $\{(\xbf_k,t_k)\}_{k\in[K]}$ \\ \hline
$\pmb{\lambda}$ & $[0,\infty)$ & search space for sparsity threshold $\hat{\lambda}$\\ \hline
$\tau$ & $(0,1]$ & $\psi$ decay tolerance \\ \hline
\end{tabular}
\caption{Hyperparameters for the WSINDy Algorithm \ref{alg:wsindypde}. Alternatively, $\mbf$ can be automatically selected from the data using the method in Appendix \ref{app:corner}. Note that the number of query points $K$ is determined from $\mbf$ and $\sbf$ using \eqref{numK}.}
\label{hypparam}
\end{table}
%

%\begin{minipage}{0.92\textwidth}
%\noindent $\what =$ \textbf{WSINDy\_PDE}$(\Ubf, \xbf,\tbf ;\ (f_j)_{j\in [J]},\, \boldsymbol{\alpha}, \abf,\sbf,\,\Sigma,\,\gamma,\,\lambda)$:
%\begin{enumerate}
%\item For $s=0, \dots, S$: compute vectors $(\phi_d^{(\alpha_d^s)}(\xbf_d))_{d\in [D+1]}$ from the support sizes $\abf = (a_1, \dots, a_{D+1})$ 
%\item For $j=1, \dots, J$: compute the $(D+1)$-dimensional array $f_j(\Ubf)$
%\item Collect convolution query points $\{(x_k,t_k)\}_{k\in[K]}$ from the computational grid $(\xbf,\tbf)$ by subsampling in each dimension with frequencies $\sbf = (s_1,s_2,\dots, s_{D+1})$ (i.e.\ for each $d\in[D+1]$ keep every $s_d$ point along the $d$th coordinate for a total of $(N_d-2a_d)/s_d$ equally-spaced points from $\xbf_d$)
%\item Compute Gram matrix $\Gbf$ and right-hand side $\bbf$ given by 
%\[\Gbf_{k,(s-1)J+j} =  \Psi^s * f_j(\Ubf)(x_k,t_k), \quad \bbf_k = \Psi^0 * \Ubf(x_k,t_k)\]
%utilizing the separable structure
%\[D^{\alpha^s}\psi(\xbf,\tbf) = \phi_1^{(\alpha_1^s)}(\xbf_1)\otimes\cdots\otimes \phi_{D+1}^{(\alpha_1^s)}(\tbf)\]
%\item Solve the generalized least-squares problem with $\ell_2$-regularization 
%\[\what = \text{argmin}_{\wbf}\left\{(\Gbf\wbf-\bbf)^T\covm^{-1}(\Gbf\wbf-\bbf)+\gamma^2\nrm{\wbf}_2^2\right\},\] 
%using sequential thresholding with parameter $\lambda$ to enforce sparsity.\\[10pt] 
%\end{enumerate}
%\end{minipage}\\

 %(i.e.\ for each $d\in[D+1]$ keep every $s_d$ point along the $d$th coordinate for a total of $(N_d-2a_d)/s_d$ equally-spaced points from $\xbf_d$)
%

\section{Examples}\label{sec:examples}
\begin{table}[H]
{\centering
\begin{tabular}{|l|l|}
\hline
\text{Inviscid Burgers }\qquad &$\partial_t u = -\frac{1}{2} \partial_x (u^2)$ \\ \hline
\text{Korteweg-de Vries (KdV) }\qquad &$\partial_t u  = -\frac{1}{2} \partial_x (u^2) -\partial_{xxx}u$ \\\hline
\text{Kuramoto-Sivashinsky (KS)}\qquad &$\partial_t u  = -\frac{1}{2} \partial_x (u^2) -\partial_{xx}u -\partial_{xxxx} u$\\\hline
\text{Nonlinear Schr\"odinger (NLS)}\qquad &$\begin{cases}\partial_t u= \frac{1}{2}\partial_{xx}v +u^2v+v^3 \\ \partial_t v = -\frac{1}{2}\partial_{xx} u-uv^2-u^3\end{cases}$\\\hline
\text{Sine-Gordon (SG)}\qquad &$\partial_{tt} u  = \vphantom{\frac11} \partial_{xx}u+\partial_{yy}u - \sin(u)$\\\hline
\text{Reaction-Diffusion (RD)}\qquad &$\begin{cases} \partial_t u  = \frac{1}{10} \partial_{xx}u+\frac{1}{10} \partial_{yy}u -uv^2-u^3+v^3+u^2v+u \\ \partial_t v  = \frac{1}{10} \partial_{xx}v+\frac{1}{10} \partial_{yy}v+v-uv^2-u^3-v^3-u^2v \end{cases}$ \\\hline
\text{2D Navier-Stokes (NS)}\qquad & $\partial_t \omega = -\partial_x(\omega u)-\partial_y(\omega v) +\frac{1}{100} \partial_{xx}\omega+\frac{1}{100} \partial_{yy}\omega$ \\\hline
\end{tabular}}
\caption{PDEs used in numerical experiments, written in the form identified by WSINDy. Note that domain specification and boundary conditions are given in Appendix \ref{app:numpde}. }
\label{pdetable}
\end{table}
We now demonstrate the effectiveness of WSINDy by recovering the seven PDEs listed in Table \ref{pdetable} over a range of noise levels, amplitudes and model libraries. The examples below show that WSINDy provides orders of magnitude improvements over derivative-based methods \cite{rudy2017data}, with reliable and accurate recovery of four out of the seven PDEs under noise levels as high as $100\%$ (defined in \eqref{sigmadef} and \eqref{snrdef}) and for all PDEs under $20\%$ noise. In contrast to the weak recovery methods in \cite{reinbold2020using,gurevich2019robust}, WSINDy uses (i) the convolutional weak form \eqref{conv_disc} and FFT-based implementation \eqref{convfft}, (ii) improved thresholding and automatic selection of the sparsity threshold $\hat{\lambda}$ via \eqref{MSTLS1} and \eqref{MSTLS2}, and (iii) rescaling using \eqref{scalext} and \eqref{scaleu}). The effects of these improvements are discussed in Sections \ref{sec:results} and \ref{sec:accresults}.

To test robustness to noise, a noise ratio $\sigma_{NR}$ is specified and a synthetic ``observed'' dataset 
\[\Ubf = \Ubf^\star+\epsilon\]
is obtained from a simulation $\Ubf^\star$ of the true PDE\footnote{Details on the numerical methods and boundary conditions used to simulate each PDE can be found in Appendix \ref{app:numpde}.} by adding white noise with variance $\sigma^2$ to each data point, where 
\begin{equation}\label{sigmadef}
\sigma :=\sigma_{NR}\nrm{\Ubf^\star}_{RMS}:=\sigma_{NR} \left(\frac{1}{(N_1\cdots N_DN_{D+1})}\sum_{k_1=1}^{N_1}\cdots\sum_{k_{D+1}=1}^{N_{D+1}}\left(\Ubf^\star_{k_1,\dots,k_{D+1}}\right)^2\right)^{1/2}.
\end{equation}
We examine noise ratios $\sigma_{NR}$ in the range $[0,1]$ and often refer to the noise level as $\sigma_{NR}$ or say that the data contains $100\sigma_{NR}\%$ noise. We note that the resulting true noise ratio
\begin{equation}\label{snrdef}
\sigma_{NR}^\star:= \frac{\nrm{\ep}_{RMS}}{\nrm{\Ubf^\star}_{RMS}}
\end{equation}
matches the specified $\sigma_{NR}$ to at least four significant digits in all cases and so we only list $\sigma_{NR}$. In the cases where the state variable itself is multi-component, as in the nonlinear Schr\"odinger equation, reaction-diffusion system, and Navier-Stokes (see Table \ref{pdetable}), a separate variance $\sigma^2$ is used to compute the noise $\epsilon$ in each component, so that $\sigma_{NR}$ is the same in each component.

\subsection{Performance Measures}

To measure the ability of the algorithm to correctly identify the terms having nonzero coefficients, we use the \textit{true positivity ratio} (introduced in \cite{LagergrenNardiniMichaelLavigneEtAl2020ProcRSocA}) defined by 
\begin{equation}\label{tpr}
\text{TPR}(\what) = \frac{\text{TP}}{\text{TP}+\text{FN}+\text{FP}}
\end{equation}
where TP is the number of correctly identified nonzero coefficients, FN is the number of coefficients falsely identified as zero, and FP is the number of coefficients falsely identified as nonzero. Identification of the true model results in a TPR of 1, while identification of half of the correct nonzero terms and no falsely identified nonzero terms results in TPR of 0.5 (e.g.\ the 2D Euler equations $\partial_t\omega = -\partial_x(\omega u) -\partial_y(\omega v)$ result in a TPR of 0.5 if the underlying true model is the 2D Navier-Stokes vorticity equation). We will see that in several cases the TPR remains above $0.95$ even as the noise level approaches $1$. The loss function $\CalL(\lambda)$ (defined in \eqref{lossfcn}) and the resulting learned sparsity threshold $\widehat{\lambda}$ (defined in \eqref{MSTLS2}) provide additional information on the algorithm's ability to identify the correct model terms with respect to the noise level. In particular, sensitivity to the sparsity threshold suggests that automatic selection of $\widehat{\lambda}$ is essential to successful recovery in the large noise regime.  

To assess the accuracy of the recovered coefficients we use two metrics. The first measures the maximum error in the true non-zero coefficients and is defined
\begin{equation}\label{err}
E_\infty(\what) := \max_{\{j\ :\ \wstar_j\neq 0\}}\frac{|\what_j - \wstar_j|}{|\wstar_j|}.
\end{equation}
$E_\infty$ determines the number of significant digits in the recovered true coefficients. We also measure the $\ell^2$ distance in parameter space using 
\begin{equation}\label{err2}
E_2(\what) := \frac{\nrm{\what-\wstar}_{RMS}}{\nrm{\wstar}_{RMS}},
\end{equation}
which provides information regarding the magnitudes of coefficients that are falsely identified as nonzero. Often when a term is falsely identified and the resulting nonzero coefficient is small, a larger sparsity factor will result in idenfitication of the true model.\\

For each system in Table \ref{pdetable} and each noise ratio $\sigma_{NR}\in\{0.025k\ :\ k\in\{0,\dots, 40\}\}$ we run WSINDy on 200 instantiations of noise\footnote{We find that 200 runs sufficiently reduces variance in the results.} and average the results of error statistics \eqref{tpr}, \eqref{err}, and \eqref{err2}. Computations were carried out on the University of Colorado Boulder Blanca Condo cluster \footnote{2X Intel Xeon 5218 at 2.3 GHz with 22 MB cache, 16 cores per cpu, and 384 GB ram.}.

\subsection{Implementation Details}

The hyperparameters used in WSINDy applied to each of the PDEs in Table \ref{pdetable} are given in Table \ref{hptable}. To select test function discrete support lengths we used a combination of manual tuning and the changepoint method\footnote{For Burgers, KdV, and KS we set $\widehat{\tau} = 3$ (defined in Appendix \ref{app:corner}) while for NLS, SG, RD and NS we used $\widehat{\tau} = 1$. For Kuramoto-Sivashinsky and nonlinear Schr\"odinger's we chose $(m_x,m_t)$ values nearby that had better performance.} described in Appendix \ref{app:corner}. Across all examples the decay tolerance for test functions is fixed at $\tau = 10^{-10}$ and the search space $\pmb{\lambda}$ for the threshold $\widehat{\lambda}$ was fixed at 
\[\pmb{\lambda} = \left\{10^{-4+j\frac{4}{49}}\ :\ j\in \{0,\dots,49\}\right\},\]
(i.e.\ $\pmb{\lambda}$ contains 50 points with $\log_{10}(\pmb{\lambda})$ equally spaced from $-4$ to $0$). We fix the subsampling frequencies $(s_x,s_t)$ to $(\frac{N_1}{50},\frac{N_2}{50})$ for PDEs in one spatial dimension and to $(\frac{N_1}{25},\frac{N_3}{25})$ for two spatial dimensions, where the dimensions $(N_1,N_2,N_3)$ depend on the dataset. Additional information about the convolutional weak discretization is included in Table \ref{spectable}, such as the dimensions and condition number of the rescaled Gram matrix $\tilde{\Gbf}$ (computed from a dataset with $20\%$ noise), test function polynomial degrees $(p_x,p_t)$, scale factors $(\gamma_u,\gamma_x,\gamma_t)$, and start-to-finish walltime of Algorithm \ref{alg:wsindypde} with all computations performed serially on a laptop with an 8-core Intel i7-2670QM CPU with 2.2 GHz and 8 GB of RAM.

{\footnotesize
\begin{table}
\centering
\begin{tabular}{|c||c|c|c|c|c|}
\hline
PDE & $\Ubf$ & $f_j$ & $\pmb{\alpha}$ & $(m_x,m_t)$ & $(s_x,s_t)$ \\ \hline \hline
Burgers & $256\times 256$ & $(u^{j-1})_{j\in[7]}$ & $((\ell,0))_{0\leq \ell \leq 6}$ & $(60,60)$ & $(5,5)$\\ \hline
KdV & $400\times 601$ & $(u^{j-1})_{j\in[7]}$ & $((\ell,0))_{0\leq \ell \leq 6}$ & $(45,80)$ & $(8,12)$ \\ \hline
KS & $256\times 301$ & $(u^{j-1})_{j\in[7]}$ & $((\ell,0))_{0\leq \ell \leq 6}$ & $(23,22)$ & $(5,6)$ \\ \hline
NLS & $2\times 256\times 251$ & $(u^nv^m)_{0\leq n+m \leq 6}$ & $((\ell,0))_{0\leq \ell \leq 6}$ & $(19,25)$ & $(5,5)$\\ \hline
SG & $129\times 403\times 205 $ &$(u^{n-1})_{n\in[5]}$, $(\sin(mu),\cos(mu))_{m=1,2}$ & $((\ell,0,0),(0,\ell,0))_{0\leq \ell \leq 4}$ & $(40,25)$ & $(5,8)$ \\ \hline
RD & $2\times 256\times 256\times 201$ & $(u^nv^m)_{0\leq n+m \leq 4}$ & $((\ell,0,0),(0,\ell,0))_{0\leq \ell \leq 5}$ & $(13,14)$ & $(13,12)$ \\ \hline
NS & $3\times 324\times 149\times 201$ & $\begin{dcases} 
(\omega^nu^mv^q)_{0\leq n+m+q \leq 2}, &  |\alpha^s| = 0 \\
(\omega^nu^mv^q)_{0\leq n+m+q \leq 3,n>0}, & |\alpha^s| > 0\end{dcases}
$ & $((\ell,0,0),(0,\ell,0))_{0\leq \ell \leq 2}$ & $(31,14)$ & $(12,8)$ \\ \hline
\end{tabular}
\caption{WSINDy hyperparameters used to identify each example PDE.}
\label{hptable}
\end{table}}

{\footnotesize
\begin{table}
\centering
\begin{tabular}{|c||c|c|c|c|c|}
\hline
PDE & $\tilde{\Gbf}$ & $\kappa(\tilde{\Gbf})$ & $(p_x,p_t)$ & $(\gamma_u,\gamma_x,\gamma_t)$ & Walltime (sec) \\ \hline \hline
Burgers & $784\times 43$ &  $9.6\times 10^4$ & $(7,7)$ & $(4.5\times 10^{-4},0.0029,1.1)$ & $0.12$ \\ \hline
KdV & $1443\times 43$ & $2.6\times 10^5$  & $(8,7)$ & $(5.7\times 10^{-4},8.3,1250)$ & $0.39$ \\ \hline
KS & $1806\times 43$ &  $9.4\times 10^3$  & $(10,10)$ & $(0.26,0.74,0.091)$ & $0.24$ \\ \hline
NLS & $1804\times 190$ & $6.0\times 10^4$ & $(11,10)$ & $(0.33,3.1,9.4)$ & $2.5$ \\ \hline
SG & $13000\times 73$ & $6.8\times 10^3$ & $(8,10)$ & $(0.23,8.1,8.1)$ & $29$\\ \hline
RD & $11638\times 181$ & $4.5\times 10^3$ & $(13,12)$& $(0.86,6.5,1.4)$ & 75 \\ \hline
NS & $3872\times 50$ &  $8.2\times 10^2$  & $(9,12)$& $(0.53,0.72,2.4)$ & 12\\ \hline
\end{tabular}
\caption{Additional specifications resulting from the choices in Table \ref{hptable}. The last column shows the start-to-finish walltime of Algorithm \ref{alg:wsindypde} with all computations in serial measured on a laptop with an 8-core Intel i7-2670QM CPU with 2.2 GHz and 8 GB of RAM.}
\label{spectable}
\end{table}}

\subsection{Comments on Chosen Examples}

The primary reason for choosing the examples in Table \ref{pdetable} is to demonstrate that WSINDy can successfully recover models over a wide range of physical phenomena such as spatiotemporal chaos, nonlinear waves, shock-forming solutions, and complex limit cycles. 

\begin{figure}
%\begin{tabular}{cc}
%		\includegraphics[trim={15 10 5 15},clip,width=0.45\textwidth]{burgers_clean} & 
%		\includegraphics[trim={15 10 5 15},clip,width=0.45\textwidth]{burgers_noise} 
%\end{tabular}
\includegraphics[trim={40 0 40 10},clip,width=0.8\textwidth]{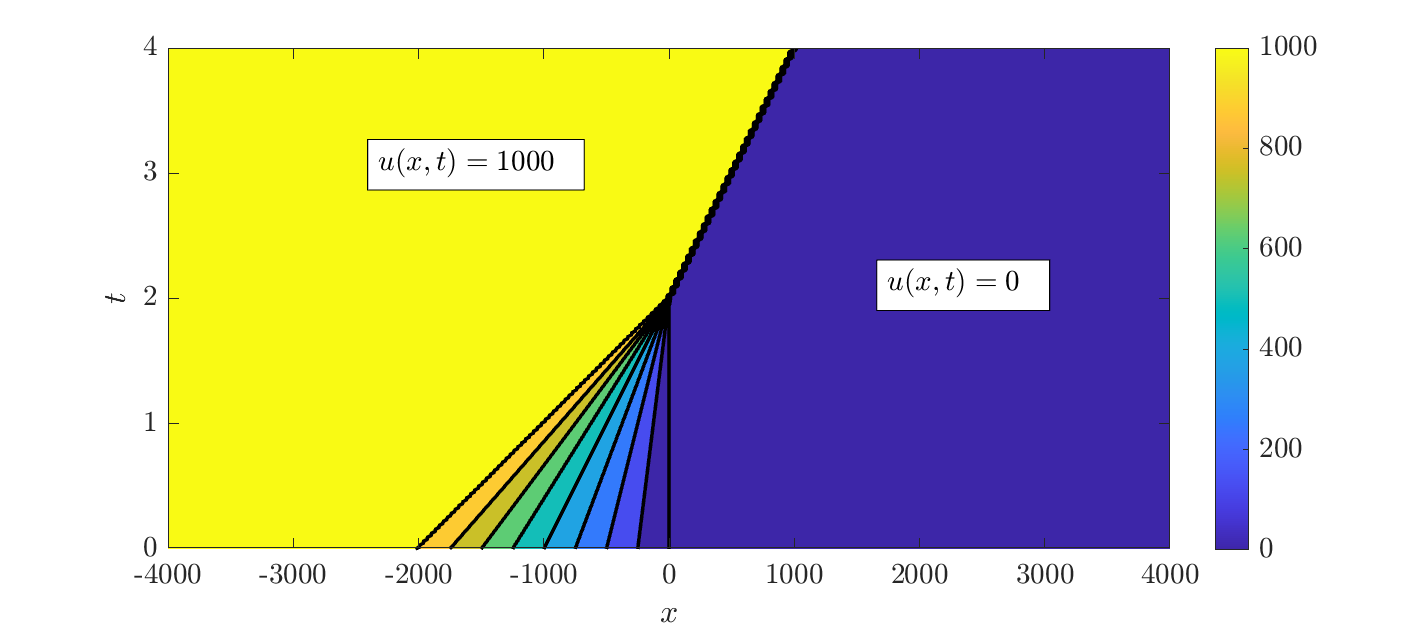}
\caption{Characteristics of the shock-forming solution \eqref{burgers} used to identify the inviscid Burgers equation. A shock forms at time $t = 2$ and travels along the line $x = 500(t-2)$.}
\label{burgerssols}
\end{figure}

Recovery of the inviscid Burgers equation demonstrates (i) that WSINDy can discover PDEs from solutions that can only be understood in a weak sense\footnote{This is conjectured in \cite{gurevich2019robust} with suggestions for how to modify the test functions in order to integrate discontinuous data but is not carried out. We demonstrate here that no such modification is necessary.} and (ii) that discovery in this case is just as accurate and robust to noise and scaling as with smooth data (i.e.\ no special modifications of the algorithm are required to discover models from discontinuous data). We use an analytical weak solution with continuous initial data that becomes discontinuous in finite time and forms a shock that propagates with constant speed (see Figure \ref{burgerssols} for plots of the characteristic curves). In addition, both the inviscid Burgers and KdV equations demonstrate that WSINDy successfully recovers the correct models for nonlinear transport data with large amplitude. Both datasets have mean amplitudes on the order of $10^3$ (in addition KdV is given over a time window of $T = 10^{-3}$), and hence are not identifiable from large polynomial libraries using naive approaches. The sparsification and rescaling measures in Sections \ref{sec:sparse} and \ref{sec:rescale} remove this barrier. 

The Sine-Gordon equation\footnote{We have not included experiments involving multiple-soliton solutions to Sine-Gordon, however the success of WSINDy applied to KdV, nonlinear Schr\"odinger and Sine-Gordon suggests that the class of integrable systems could be a fruitful avenue for future research.} is used to show both that trigonometric library terms can easily be identified alongside polynomials and that hyperbolic problems do not seem to present further challenges. Discovery of the Sine-Gordon equation also appears to be particularly robust to noise, which suggests that the added complexity of having multiple spatial dimensions is not in general a barrier to identification. 

For the nonlinear Schr\"odinger and reaction-diffusion systems, we test the ability of WSINDy to select the correct monomial nonlinearities from an excessively large model library. Using a library of 190 terms for nonlinear Schr\"odinger's and 181 terms for reaction-diffusion (see the dimensions of $\tilde{\Gbf}$ in Table \ref{spectable}), we demonstrate successful identification of the correct nonzero terms. Moreover, for the reaction-diffusion system misidentified terms directly reflect the existence of a limit cycle\footnote{We note that discovery of the same reaction-diffusion system from a much smaller library of terms is shown in \cite{rudy2017data,reinbold2020using}, but with different initial conditions that result in a spiral wave limit cycle. Our choice of initial conditions is motivated below in Appendix \ref{app:numpde}.}.

\subsection{Results: Model Identification}\label{sec:results}

\begin{figure}
\begin{tabular}{cc}
		\includegraphics[trim={0 0 35 15},clip,width=0.48\textwidth]{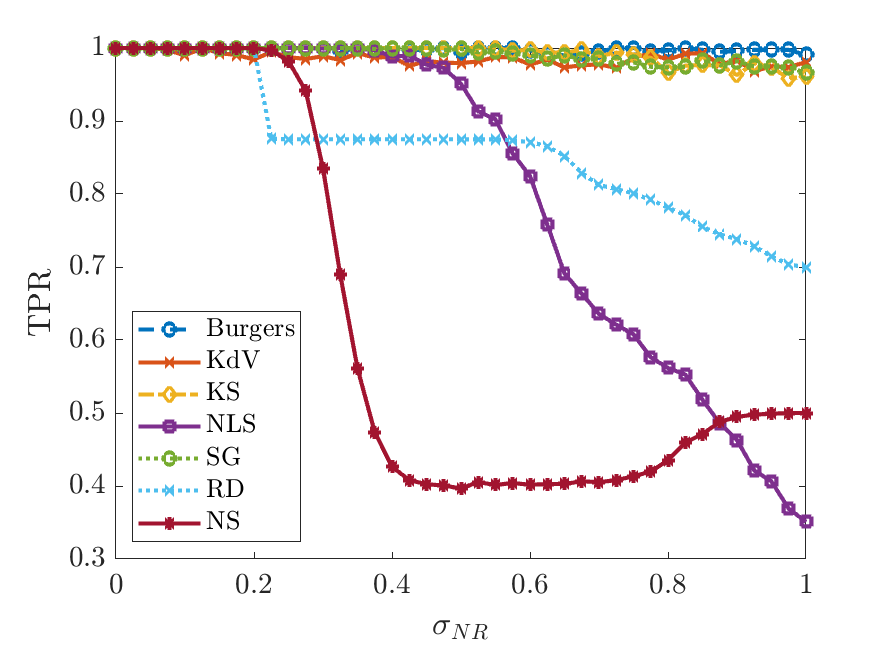} & 
		\includegraphics[trim={0 0 35 15},clip,width=0.48\textwidth]{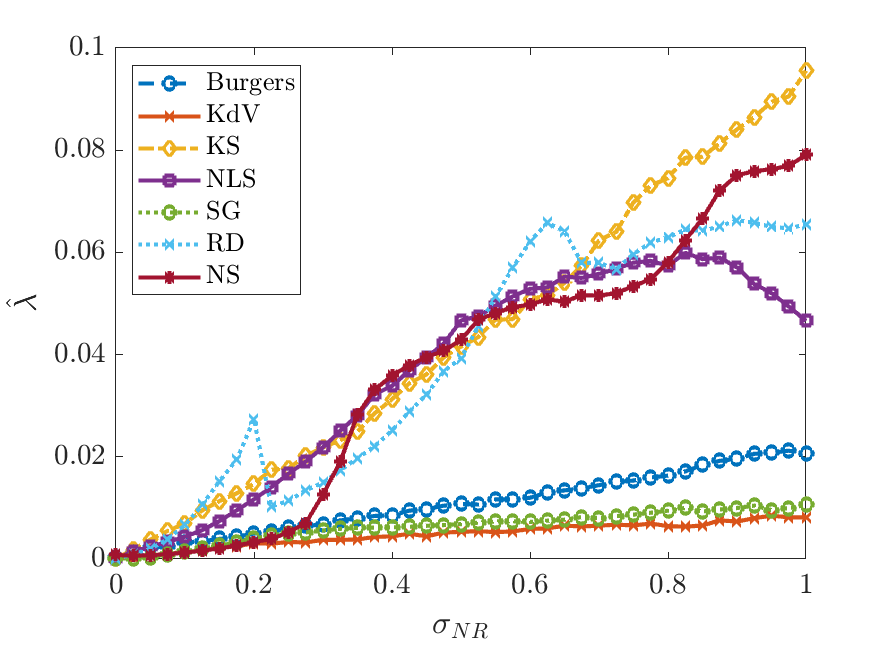}
\end{tabular}
\caption{Left: average TPR (total positivity ratio, defined in \eqref{tpr}) for each of the PDEs in Table \ref{pdetable} computed from 200 instantiations of noise for each noise level $\sigma_{NR}$. Right: average learned threshold $\widehat{\lambda}$ (defined in \eqref{MSTLS2}).}
\label{tprfig}
\end{figure}

\begin{figure}
     \centering
     \begin{subfigure}[b]{0.48\textwidth}
         \centering
         \includegraphics[trim={0 0 35 15},clip,width=\textwidth]{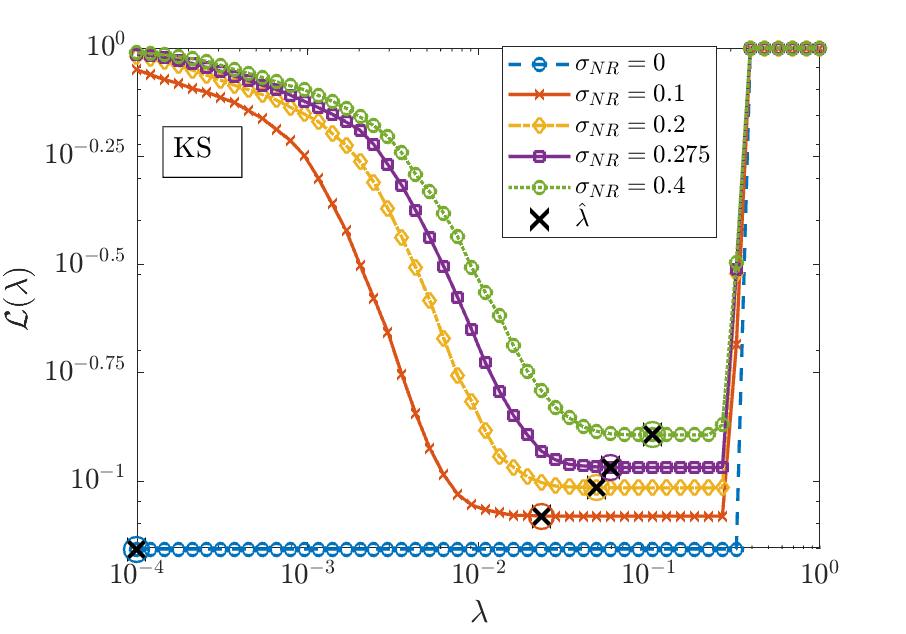}
         \caption{Kuramoto-Sivashinsky.}
         \label{fig:lossKS}
     \end{subfigure}
     \begin{subfigure}[b]{0.48\textwidth}
         \centering
         \includegraphics[trim={0 0 35 15},clip,width=\textwidth]{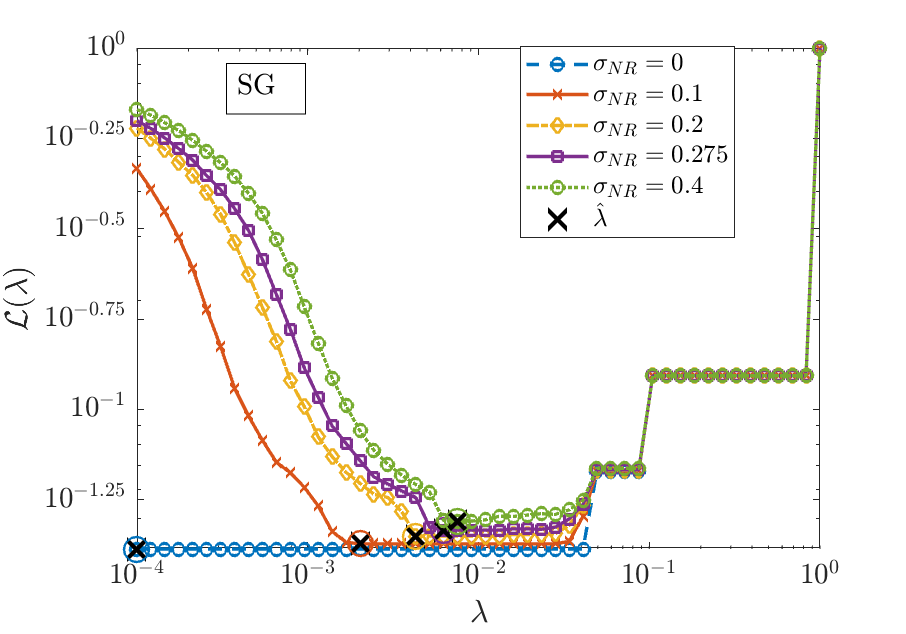}
         \caption{Sine-Gordon.}
         \label{fig:lossSG}
     \end{subfigure}
     \begin{subfigure}[b]{0.48\textwidth}
         \centering
         \includegraphics[trim={0 0 35 15},clip,width=\textwidth]{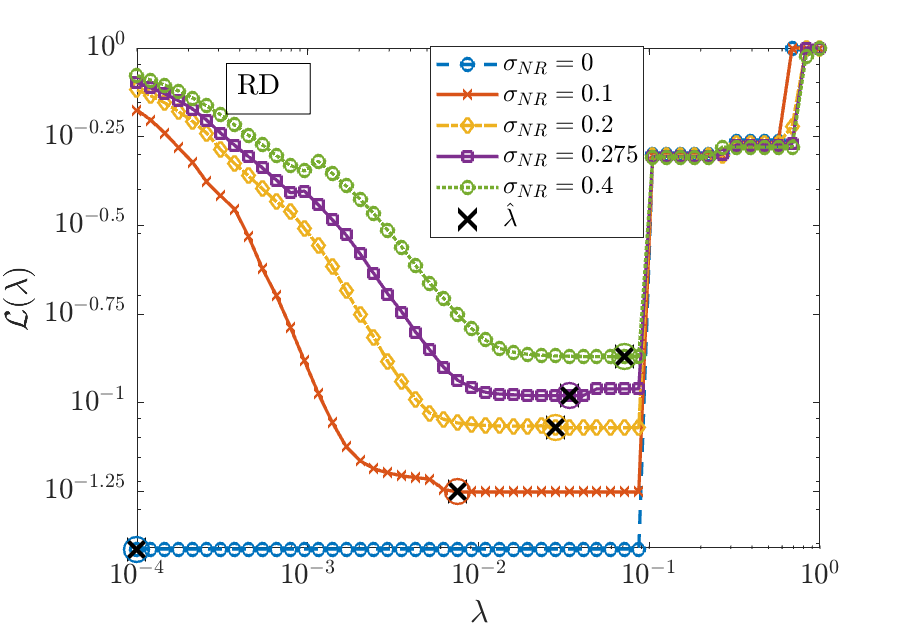}
         \caption{Reaction-diffusion.}
         \label{fig:lossRD}
     \end{subfigure}
     \begin{subfigure}[b]{0.48\textwidth}
         \centering
         \includegraphics[trim={0 0 35 15},clip,width=\textwidth]{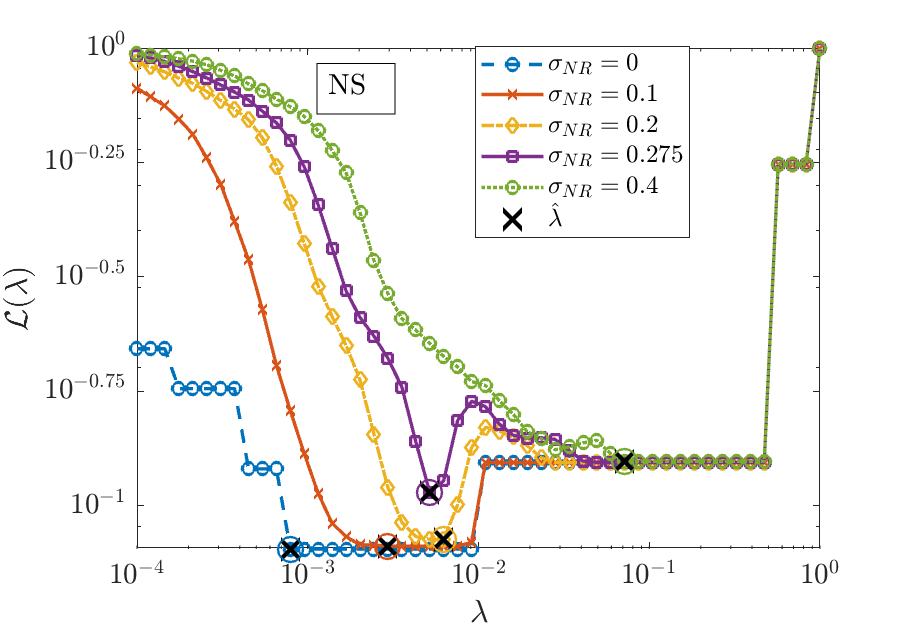}
         \caption{Navier-Stokes.}
         \label{fig:lossNS}
     \end{subfigure}
        \caption{Plots of the average loss function $\CalL(\lambda)$ and resulting optimal threshold $\widehat{\lambda}$ for the Kuramoto Sivashinsky, Sine-Gordon, Reaction diffusion and Navier-Stokes equations.}
        \label{fig:three graphs}
\label{loss}
\end{figure}

Performance regarding the identification of correct nonzero terms in each model is reported in Figures \ref{tprfig} and \ref{loss}, which include plots of the average TPR, the learned threshold $\widehat{\lambda}$, and the loss function $\CalL(\lambda)$ (defined in \eqref{tpr}, \eqref{MSTLS2}, and \eqref{lossfcn}, respectively). As we will discuss, significant decreases in average TPR are often accompanied by transitions in the identified $\widehat{\lambda}$.  

Figure \ref{tprfig} (left) shows that for inviscid Burgers, Korteweg-de Vries, Kuramoto-Sivashinsky and Sine-Gordon, the average TPR stays above 0.95 even for noise levels as high as $100\%$ (i.e.\ WSINDy reliably identifies these models in the presence of noise that has an $L^2$-norm comparable to that of the underlying clean data). The average TPR for nonlinear Schr\"odinger's stays above 0.95 until 50\% noise, a drastic improvement over previous studies \cite{rudy2017data}, after which spurious higher-degree monomials are selected. This is to be expected from the large library of 190 terms used here\footnote{Recovery of nonlinear Schr\"odinger's is significantly more robust with a smaller library.}. 

 We observe in Figure \ref{tprfig} (right) that the learned threshold $\widehat{\lambda}$ increases with $\sigma_{NR}$, suggesting that automatic selection of $\widehat{\lambda}$ in the learning algorithm \eqref{MSTLS2} is crucial to the algorithm's robustness to noise. For the Kuramoto-Sivashinsky equation in particular, which has a minimum nonzero coefficient of 0.5 (multiplying $\partial_x(u^2)$), we find that $\widehat{\lambda}$ approaches $0.1$ as $\sigma_{NR}$ approaches 1, which implies that at higher noise levels the range of $\widehat{\lambda}$ values that is necessary\footnote{By definition \eqref{MSTLS2}, $\widehat{\lambda}$ is the minimum value in $\pmb{\lambda}$ that minimizes the loss $\CalL$ \eqref{MSTLS2}, hence values in $\pmb{\lambda}$ below $\widehat{\lambda}$ are precisely the thresholds that result in misidentification of the correct model by overfitting, while thresholds above $\min_{\{j\,:\,\wstar_j\neq 0\}}|\wstar_j|$ necessarily underfit the model.} for correct model identification is approximately $(\sim 0.1,\,\sim 0.5)$. Since it is highly unlikely that this range of admissible values would be known \textit{a priori}, the chances of manually selecting a feasible $\widehat{\lambda}$ for Kuramoto-Sivashinsky are prohibitively low in the large noise regime (see Figure \ref{fig:lossKS} for visualizations of the loss $\CalL$ applied to KS data). Automatic selection of $\widehat{\lambda}$ thus removes this sensitivity. In contrast, $\widehat{\lambda}$ is largely unaffected by increases in $\sigma_{NR}$ for Burgers, Korteweg-de Vries and Sine-Gordon. In particular, Figure \ref{fig:lossSG} shows little qualitative changes in the loss landscape for Sine-Gordon in the range $0.1\leq \sigma_{NR}\leq 0.4$. 

For reaction-diffusion, the average TPR falls below 0.95 at $22\%$ noise, after which WSINDy falsely identifies linear terms in $u$ and $v$. The underlying solution settles into a limit cycle, which means that at every point in space the solution will oscillate. If the true model is given by the compact form $\partial_t\ubf = \CalA(\ubf)$ for $\ubf = [u\ v]^T$, then the misidentified model in all trials for noise levels in the range $0.25\leq \sigma_{NR}\leq 0.55$ is given by
\begin{equation}\label{oscmod}
\partial_t \ubf = \beta \CalA(\ubf)+\alpha\begin{pmatrix} 0 & 1\\ -1 & 0\end{pmatrix}\ubf
\end{equation}
for some $\alpha>0$ and $\beta\approx 1$ dependent on $\sigma_{NR}$. The falsely identified nonzero terms convey that at each point in space the solution is oscillating at a uniform frequency, but with variable amplitude and phase determined by the initial conditions\footnote{This is discussed further in Appendix \ref{sec:RD}.}. Hence, in the presence of certain lower-dimensional structures (in this case a limit cycle), higher noise levels result in a mixture of the true model with a spatially-averaged reduced model. This shift between detection of the correct model and the oscillatory version \eqref{oscmod} is also detectable in the learned threshold $\widehat{\lambda}$, which decreases at $\sigma_{NR} = 0.22$ (see RD data in Figure \ref{tprfig} (right)), and in the loss function $\CalL$ (Figure \ref{fig:lossRD}). At $\sigma_{NR}=0.275$ we see that $\CalL$ in Figure \ref{fig:lossRD} is minimized for $\lambda$ in the approximate range $(\sim 0.02,\, \sim 0.05)$ but also has a near-minimum for $\lambda \in (\sim 0.05,\,\sim 0.1)$. These two regions correspond to discovery of the oscillatory model \eqref{oscmod} and the true model, respectively, but since the true model has a slightly higher loss, model \eqref{oscmod} is selected. For $\sigma_{NR} \geq 0.4$ there is no longer (on average) a region of $\lambda$ that results in discovery of the true model, and WINSDy returns \eqref{oscmod} to compensate for noise.

For Navier-Stokes we see an averaging effect at higher noise, similar to the reaction-diffusion system. TPR drops below 0.95 for noise levels above $27\%$ with the resulting misidentified model being simply Euler's equations in vorticity form:
\[\partial_t\omega = -\partial_x(\omega u) -\partial_y(\omega v).\]
This is due primarily to the small viscosity $1/Re = 0.01$ resulting from Reynolds number $Re=100$, which prevents identification of the viscous forces at higher noise levels. Viewed from the perspective of the loss function $\CalL$, Figure \ref{fig:lossNS} shows that for noise levels under $\sigma_{NR} = 0.275$, minimizers of $\CalL$ are below $0.01$, while for higher noise levels, minimizers are above 0.01, rendering terms in the model with coefficient less than $0.01$ unidentifiable. Another impediment to discovery of Navier-Stokes is the low-accuracy simulation used for the clean dataset: in the noise-free setting, Table \ref{cleandata} shows that WSINDy recovers the model coefficients of Navier-Stokes to less than 3 significant digits in the absence of noise, which is the same level of accuracy exhibited by the method on each of the other systems with $5\%$ noise  (see Figure \ref{errs}). Nevertheless, with reliable recovery up to $27\%$ noise, WSINDy makes notable improvements on previous results (\cite{rudy2017data}).

%\begin{figure}
%\begin{tabular}{cc}
%	\includegraphics[trim={0 0 35 15},clip,width=0.48\textwidth]{lossKS.png} &
%	\includegraphics[trim={0 0 35 15},clip,width=0.48\textwidth]{lossSG.png}\\
%	\includegraphics[trim={0 0 35 15},clip,width=0.48\textwidth]{lossRD.png} &
%	\includegraphics[trim={0 0 35 15},clip,width=0.48\textwidth]{lossNS.png}
%	\end{tabular}
%\caption{Plots of the loss function $\CalL$ (defined in \eqref{lossfcn}) for reaction-diffusion (left) and Navier-Stokes (right) at various noise levels. For reaction-diffusion, the true model competes with a spurious oscillatory version, ultimately losing the status if minimizer for $\sigma_{NR}\geq 0.3$. For Navier-Stokes, minimizers of $\CalL$, initiallly below $0.01$, shift to abover $0.01$ for $\sigma_{NR}>0.3$ increases beyond $0.3$ the minimizer shifts minimizers is above $0.01$, rendering the terms $0.01\partial_{xx}\omega+0.01\partial_{yy}\omega$ unidentifiable and WINDy\_PDE instead returns the inviscid 2D Euler equations. Right: plots of the learned threshold $\widehat{\lambda}$ (defined in \eqref{MSTLS2}), which increases with $\sigma_{NR}$, indicating that sparsity must be enforced more strictly to observe the correct model at higher noise levels.}
%\label{loss}
%\end{figure}

\begin{table}
\begin{tabular}{|c|c|c|c|c|c|c|c|}
\hline
& Inviscid Burgers & KdV & KS & NLS & SG & RD & NS\\ \hline
$E_\infty$ & $4.3\times 10^{-5}$ & $3.1\times 10^{-7}$ & $8.1\times 10^{-7}$ & $9.4\times 10^{-8}$ & $4.3\times 10^{-5}$ & $3.9\times 10^{-10}$ & $1.1\times 10^{-3}$\\\hline
\end{tabular}
\caption{Accuracy of WSINDy applied to noise-free data $(\sigma_{NR}=0)$.}
\label{cleandata}
\end{table}

\subsection{Results: Coefficient Accuracy}\label{sec:accresults}

\begin{figure}
\begin{tabular}{cc}
		\includegraphics[trim={0 0 35 20},clip,width=0.48\textwidth]{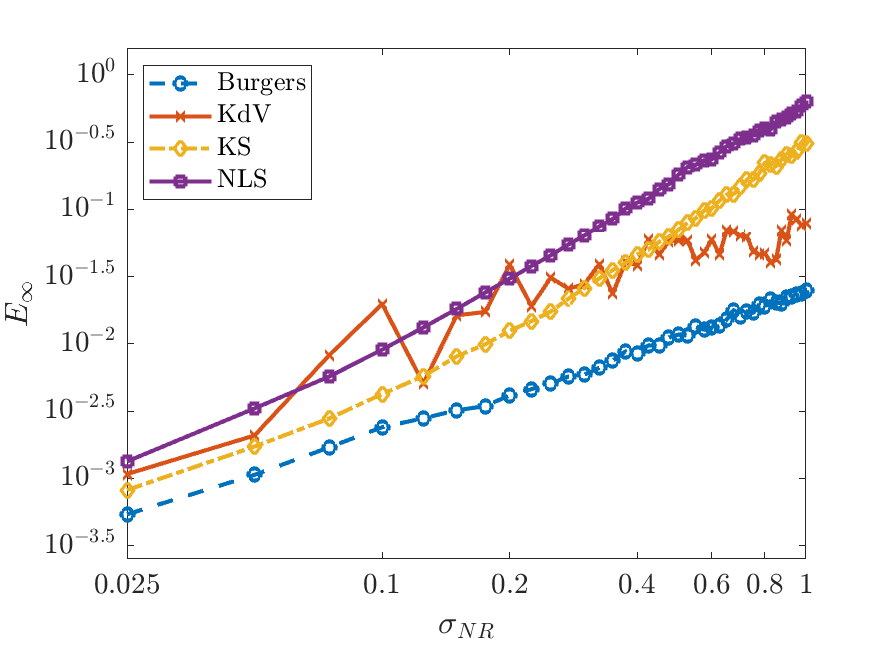} & 
		\includegraphics[trim={0 0 35 20},clip,width=0.48\textwidth]{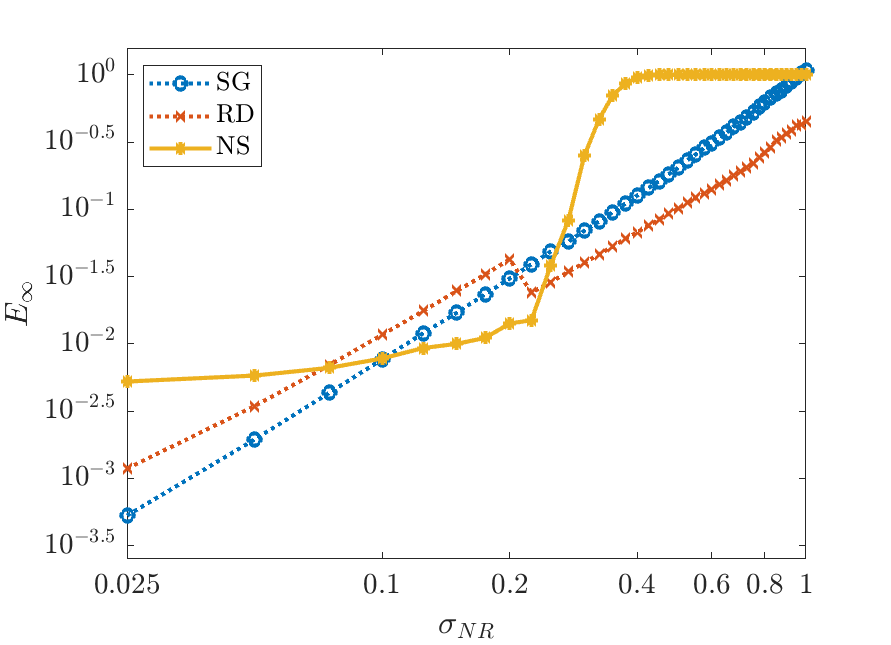} \\
		\includegraphics[trim={0 0 35 20},clip,width=0.48\textwidth]{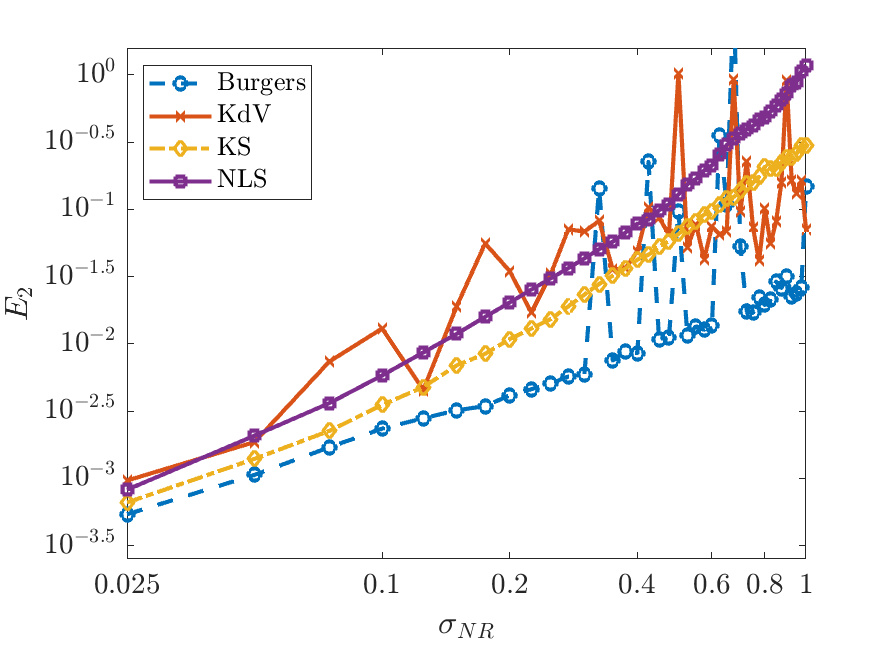} & 
		\includegraphics[trim={0 0 35 20},clip,width=0.48\textwidth]{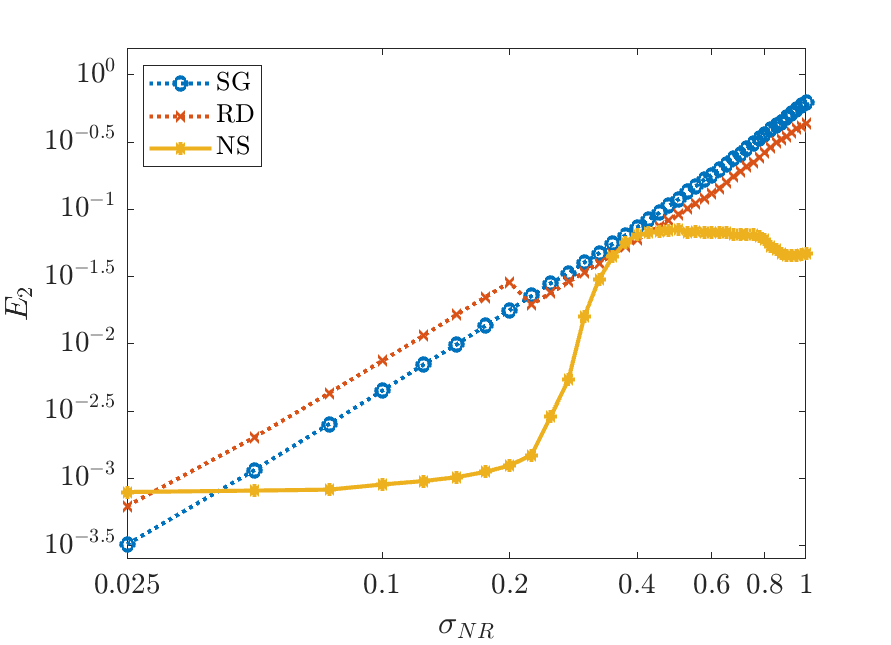} 
\end{tabular}
\caption{Coefficient errors $E_\infty$ and $E_2$ (equations \eqref{err} and \eqref{err2}) for each of the seven models Table \ref{pdetable}. Models in one and two spatial dimensions are shown on the left and right, respectively.}
\label{errs}
\end{figure}

Accuracy in the recovered coefficients is measured by $E_\infty$ and $E_2$ (defined in \eqref{err} and \eqref{err2}, respectively) and shown in Table \ref{cleandata} for $\sigma_{NR}=0$ and in Figure \ref{errs} for $\sigma_{NR}>0$. As in the ODE case, the coefficient error $E_\infty$ for smooth, noise-free data is determined by the order of accuracy of the numerical simulation method\footnote{For example, Sine-Gordon and Navier-Stokes are both integrated in time using second-order methods, hence have lower accuracy than the other examples (see Appendix \ref{app:numpde} for more details).}, since the error resulting from the trapezoidal rule is of lower order for the values $(p_x,p_t)$ used in Table \ref{spectable} (see \cite{messenger2020weak}, Lemma 1). Table \ref{cleandata} also shows that the algorithm returns reasonable accuracy for non-smooth data, with $E_\infty = 4.3\times 10^{-5}$ for the inviscid Burgers equation. 

For $\sigma_{NR}>0$,  in Figure \ref{errs} it is apparent that $E_\infty$ scales approximately as a power law $E_\infty\sim \sigma_{NR}^r$ for some $r$ approximately in the range $(\sim 1,\,\sim 2)$ in all systems except Navier-Stokes. It was observed in \cite{gurevich2019robust} that $E_\infty$ will approximately scale linearly with $\sigma_{NR}$ for Kuramoto-Sivashinsky, however our results show that in general, for larger $\sigma_{NR}$, the rate will be superlinear and dependent on the reference test function and the nonlinearities present. A simple explanation for this in the case of normally-distributed noise is the following: linear terms $\Psi^{s} * \Ubf$ will be normally-distributed with mean $\Psi^s*\Ubf^\star$ and approximate variance $\Delta x^D\Delta t\nrm{D^{\alpha^s}\psi}_2^2\sigma^2$, hence are \textit{unbiased}\footnote{In other words, equal to the noise-free case in expectation (recall that $\Ubf^\star$ is the underlying noise-free data).} and lead to perturbations that scale linearly with $\sigma_{NR}$. On the other hand, general monomial nonlinearities\footnote{With the exception of $j=2$ and odd $|\alpha^s|$, due to the fact that $\Ebb[\Psi^s*\ep^2]\approx \Ebb[\ep^2]\int_{\Omega_R} D^{\alpha^s}\psi\,dxdt= 0$.}  $\Psi^s*\Ubf^j$ with $j>1$ are {\it biased} and have approximate variance $\Delta x^D\Delta t\nrm{D^{\alpha^s}\psi}_2^2p_{2j}(\sigma)$ for $p_{2j}$ a polynomial of degree $2j$. Hence, nonlinear terms $\Psi^s*f_j(\Ubf)$ lead to biased columns of the Gram matrix $\Gbf$ with variance scaling with $\sigma^{2r}$ for some $r>1$ and proportional to $\nrm{D^{\alpha^s}\psi}_2$. Thus, for larger noise and higher-degree monomial nonlinearities, we expect superlinear growth of the error, as observed in particular with nonlinear Schr\"odinger's, Sine-Gordon, and reaction-diffusion. Nevertheless, Figure \ref{errs} suggests that a conservative estimate on the coefficient error is $E_\infty\leq \frac{\sigma_{NR}}{10}$, indicating $1-\log_{10}(\sigma_{NR})$ significant digits (e.g.\ for $\sigma_{NR} = 0.1$ we have $E_\infty\leq 10^{-2}$ for each system except KdV, indicating two significant digits), which is consistent with the ODE case \cite{messenger2020weak}.

For Burgers and Korteweg-De Vries, the average error $E_2$ at higher noise levels is affected by outliers containing a falsely-identified advection term $\partial_x u$. Since the closest pure-advection model to each of these datasets\footnote{This is found by projecting the left-hand side $\bbf$ onto the column $\partial_x \psi*\Ubf^\star$ (i.e.\ in the noise-free case).} is given by
\[\txt{0.2}{(Burgers)} \partial_t u = -(498)\partial_x u, \qquad \txt{0.2}{(KdV)} \partial_t u = -(512)\partial_x u,\]
a falsely identified $\partial_x u$ term generally has a large coefficient, whereas the true model coefficients all have magnitude 0.5 or 1. In all other cases, the values of $E_2$ and $E_\infty$ are comparable, which implies that misidentified terms do not have large coefficients and might be removed with a larger threshold. Lastly, the sigmoidal shape of $E_\infty$ and $E_2$ for Navier-Stokes is due again to the non-identification of diffusive terms at larger noise. It is interesting to note that for $\sigma_{NR}\leq 0.27$ the coefficient error for Navier-Stokes is relatively constant, in contrast to the other systems, and does not exhibit a power-law. However, at present, we do not have a concrete explanation for this behavior.

\section{Conclusion}\label{sec:conclusion}

We have extended the WSINDy algorithm to the setting of PDEs for the purpose of discovering models for spatiotemporal dynamics without relying on pointwise derivative approximations, black-box closure models (e.g.\ deep neural networks), dimensionality reduction, or other noise filtering. We have provided methods for learning many of the algorithm's hyperparameters directly from the given dataset, and in the case of the threshold $\widehat{\lambda}$, demonstrated the necessity of avoiding manual hyperparameter tuning. The underlying convolutional weak form \eqref{conv_form} allows for efficient implementation using the FFT. This naturally leads to a selection criterion for admissable test functions based on spectral decay, which is implemented in the examples above. In addition, we have shown that by utilizing scale invariance of the PDE together with a modified sparsification measure, models may be recovered from large candidate model libraries and from data that is poorly-scaled. When unsuccessful, WSINDy appears to discover a nearby sparse model that captures the dominant spatiotemporal behavior (see the discussions surrounding misidentification of the reaction-diffusion and Navier-Stokes equations in Section \ref{sec:results}). 

We close with a summary of possible future directions. In Section \ref{sec:piecewisepoly} we discussed the significance of decay properties of test functions in real and in Fourier space, as well as general test function regularity. We do not make any claim that the class $\CalS$ defined by \eqref{testfcn} is optimal, but it does appear to work very well, as demonstrated above (as well as in the ODE setting \cite{messenger2020weak}) and also observed in \cite{reinbold2020using,gurevich2019robust}. A valuable tool for future development of weak identification schemes would be the \textit{identification of optimal test functions}. A preliminary step in this direction is our use of the changepoint method described in Appendix \ref{app:corner}.

In the ODE setting, adaptive placement of test functions provided increased robustness to noise. Convolution query points can similary be strategically placed near regions of the dynamics with high information content, which may be crucial for model selection in higher dimensions. Defining regions of high information content and \textit{adaptively placing query points} accordingly would allow for identification from smaller datasets.
 
Ordinary least squares makes the assumption of i.i.d.\ residuals and should be replaced with generalized least squares to accurately reflect the true error structure. The current framework could be vastly improved by incorporating more precise statistical information about the linear system $(\Gbf,\bbf)$. The first step in this direction is the derivation of an \textit{approximate covariance matrix} as in WSINDy for ODEs \cite{messenger2020weak}.

 Accuracy in the recovered coefficients is still not entirely understood and is needed to derive recovery guarantees. It is claimed in \cite{gurevich2019robust} that at higher noise levels the scaling will approximately be linear in $\sigma_{NR}$, while we have demonstrated that this is not the case in general: the scaling depends on the nonlinearities present in the true model, the decay properties of the test functions, and accuracy of the underlying clean data. Analysis of \textit{coefficient error dependence} (on noise, amplitudes, number of datapoints, etc.) could occur in tandem with development of a generalized least-squares framework.

The examples above show that WSINDy is very robust to noise for problems involving nonlinear waves (Burgers, Korteweg de-Vries, nonlinear Schr\"odinger, Sine-Gordon) and spatiotemporal chaos (Kuramoto-Sivashinsky), but less so for data with limit cycles (reaction-diffusion, Navier-Stokes). Further, identification of Burgers, Korteweg de-Vries, and Sine-Gordon appears robust to changes in the sparsity threshold $\widehat{\lambda}$ (see Figure \ref{tprfig} (right)). A \textit{structural identifiability criteria} for measuring uncertainty in the recovery process based on identified structures (transport processes, mixing, limit cycles, etc.) would be invaluable for general model selection.

\section{Acknowledgements}
This research was supported in part by the NSF/NIH Joint DMS/NIGMS Mathematical Biology Initiative grant R01GM126559 and in part by the NSF Computing and Communications Foundations Division grant CCF-1815983. This work also utilized resources from the University of Colorado Boulder Research Computing Group, which is supported by the National Science Foundation (awards ACI-1532235 and ACI-1532236), the University of Colorado Boulder, and Colorado State University. Code used in this manuscript is publicly available on GitHub at \url{https://github.com/dm973/WSINDy_PDE}. The authors would like to thank Samuel Rudy, Kadierdan Kaheman, and Zofia Stanley for helpful discussion.

\newpage
\bibliographystyle{plain}
\bibliography{researchCU}

\appendix

\section{Learning Test Functions From Data}\label{app:corner}

\begin{figure}
\begin{tabular}{cc}
\includegraphics[trim={15 5 25 15},clip,width=0.48\textwidth]{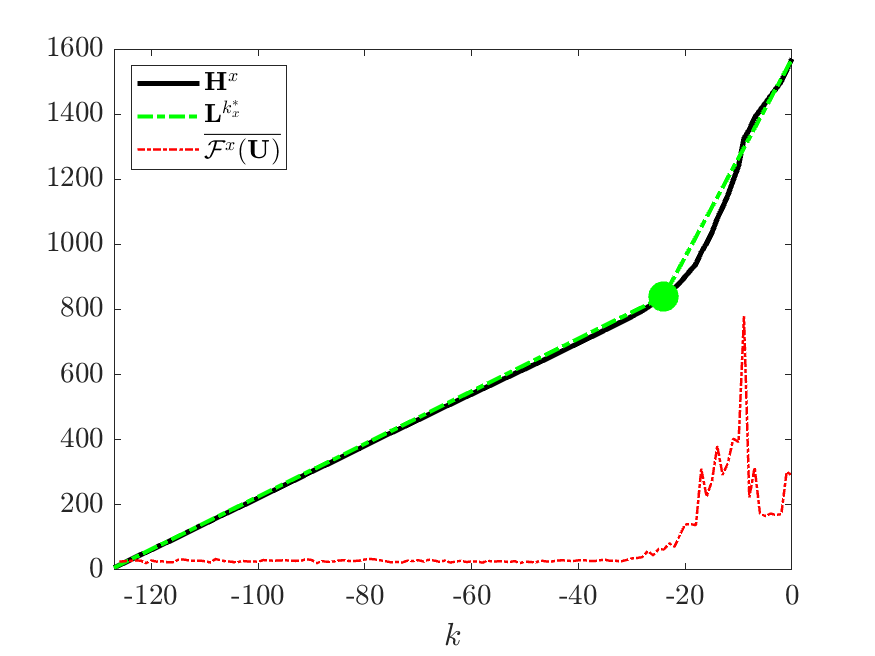} & 
\includegraphics[trim={25 10 30 20},clip,width=0.48\textwidth]{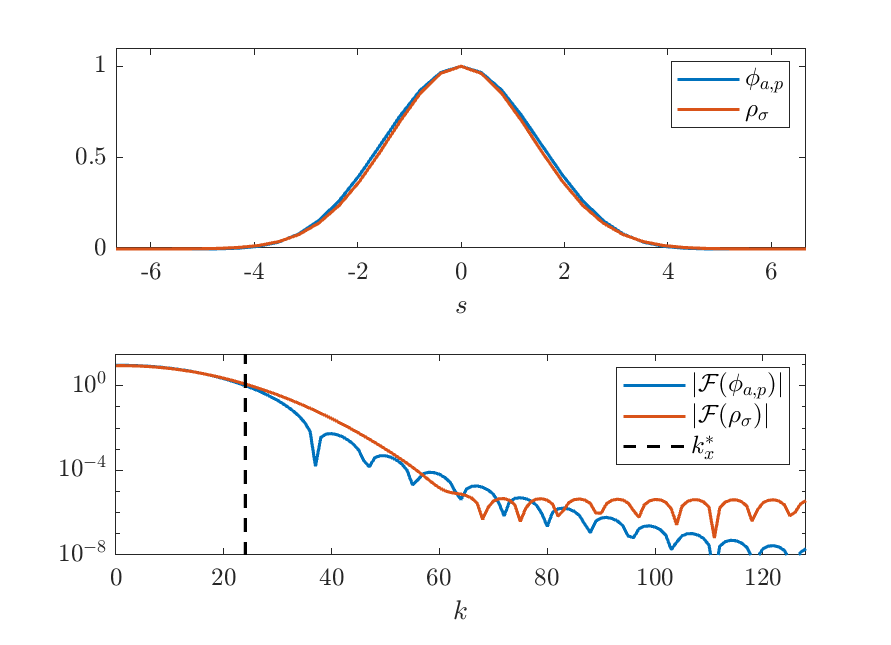}
\end{tabular}
\caption{Visualization of the changepoint algorithm for KS data with $50\%$ noise. Left: $\Hbf^x$ (defined in \eqref{cp1}) and best two-piece approximation $\Lbf^{k^*_x}$ along with resulting changepoint $k^*_x= 24$ marked in green. The noise-dominated region of $\Hbf^x$ ($k< -24$) is approximately linear as expected from the i.i.d. noise. (The time-averaged power spectrum $\overline{|\CalF^x(\Ubf)|}$ is shown in red and magnified for scale). Right: resulting test function $\phi_x = \phi_{a,p}$ and power spectrum $|\CalF(\phi_{a,p})|$ along with reference Gaussian $\rho_{\sigma}$ with $\sigma = a/\sqrt{2p+3}$. The power spectra $|\CalF(\phi_{a,p})|$ and $|\CalF(\rho_\sigma)|$ are in agreement over the signal-dominated modes $(k\leq 24)$.  (Note that the power spectrum is symmetric about zero.)}
\label{corneralg}
\end{figure}

\begin{figure}
\centering 
\includegraphics[trim={60 40 55 40},clip,width=0.8\textwidth]{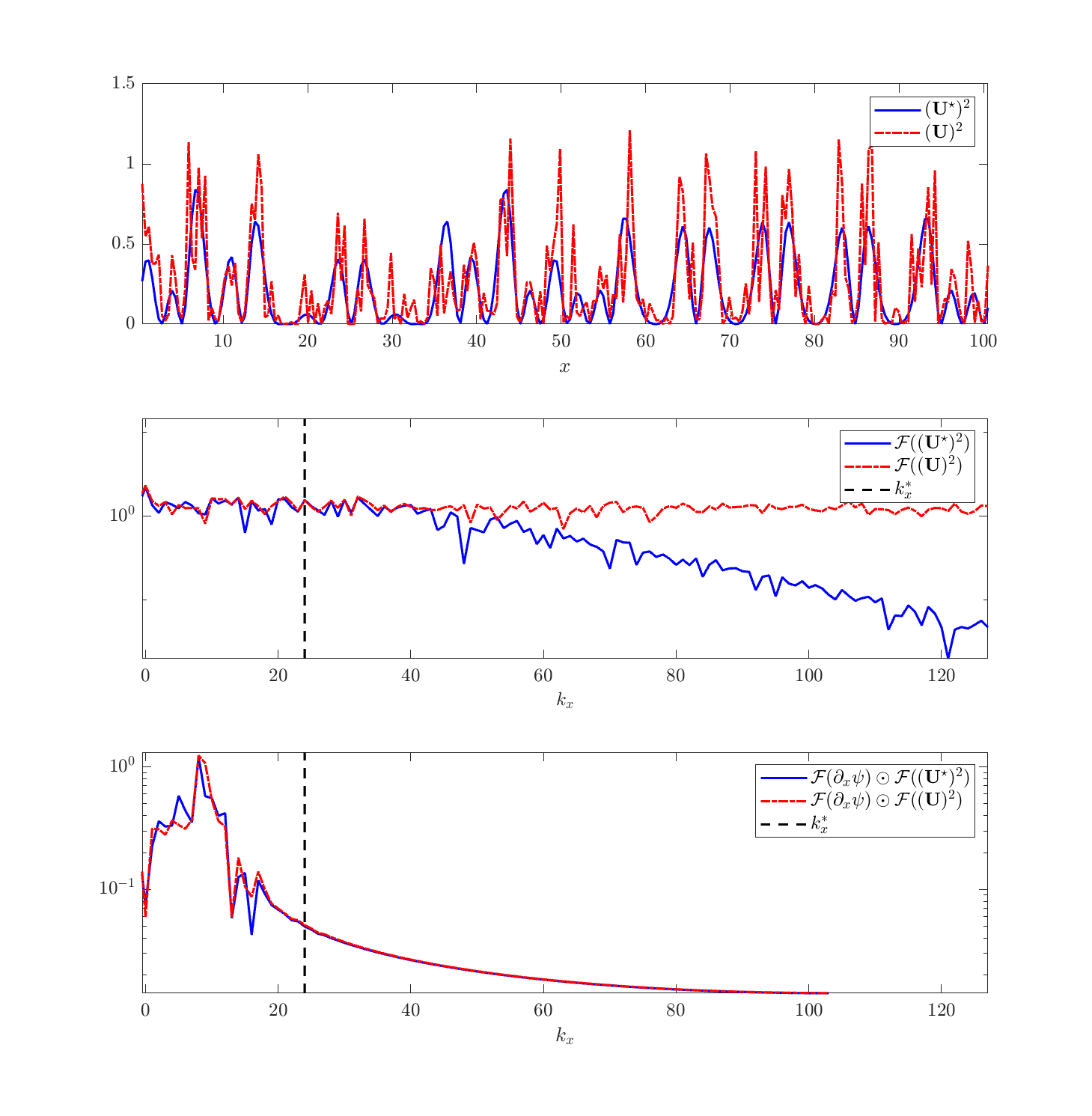} \\ 
\caption{Illustration of the test function learning algorithm using computation of $\partial_x\psi*(\Ubf^2)$ along a slice in $x$ at fixed time $t=99$ for the same dataset used in Figure \ref{corneralg}. From top to bottom: (i) clean vs. noisy data, (ii) power spectra of the clean vs. noisy data along with the learned corner point $k^*_x$, (iii) power spectra of the element-wise products $\CalF(\partial_x\psi)\odot \CalF((\Ubf^\star)^2)$ and $\CalF(\partial_x\psi)\odot \CalF((\Ubf)^2)$ (recall that these computations are embedded in the FFT-based convolution \eqref{convfft}).}
\label{visKS}
\end{figure}

Automatic selection of test functions involves two steps: (1) estimation of critical wavenumbers $(k_1^*,k_2^*,\dots)$ separating noise- and signal-dominated modes in each coordinate and (2) enforcing decay in real and in Fourier space. We will describe the process for detecting $k_x^*=k_1^*$ on data $\Ubf\in \Rbb^{N_1\times N_2}$ given over the one-dimensional spatial grid $\xbf\in \Rbb^{N_1}$ at timepoints $\tbf\in\Rbb^{N_2}$. Figures \ref{corneralg}-\ref{visKS} then illustrate the algorithm using Kuramoto-Sivashinsky data with $50\%$ noise. Below $\CalF^x$ and $\CalF^t$ denote the discete Fourier transform along the $x$ and $t$ coordinates, while $\CalF$ denotes the full two-dimensional discrete Fourier transform.

\textit{1. Detection of Critical Wavenumbers.} Assume the data has additive white noise $\Ubf = \Ubf^\star+ \ep$ with $\ep \sim \CalN(0,\sigma^2)$ and that $\CalF(\Ubf^\star)$ decays. The power spectrum of the noise $|\CalF^x(\ep)|$ is then i.i.d, hence as discussed in Section \ref{sec:piecewisepoly}, there will be a critical wavenumber $k^*_x$ in the power spectrum of the data $\CalF^x(\Ubf)$ after which the modes become noise-dominated. To detect $k_x^*$, we collapse $|\CalF^x(\Ubf)|$ into a one-dimensional array by averaging in time and then take the cumulative sum in $x$:
\begin{equation}\label{cp1}
\Hbf^x_k := \sum_{j=-N_1/2}^k \overline{|\CalF^x_j(\Ubf)|}
\end{equation}
where $\overline{|\CalF^x_j(\Ubf)|}$ is the time-average of the $j$th mode of the discrete Fourier transform along the $x$-coordinate. Since $|\CalF^x(\ep)|$ is i.i.d., $\Hbf^x$ will be approximately linear over the noise-dominated modes, which is an optimal setting for locating $k_x^*$ as a \textit{changepoint}, or in other words the corner point of the best piecewise-linear approximation\footnote{In the weighted least-squares sense with weights $\omega_k = |\Hbf^x_k|^{-1}$.} to $\Hbf^x$ using two pieces (see Figure \ref{corneralg}). An algorithm for this is given in \cite{killick2012optimal} and implemented in MATLAB using the function \texttt{findchangepts}.

\textit{2. Enforcing Decay.} We find hyperparameters for the coordinate test functions $\phi_x$ and $\phi_t$ by enforcing that (i) the changepoints $(k^*_x,k^*_t)$ are approximately $\widehat{\tau}$ standard deviations into the tail of the spectra $\widehat{\phi}_x$ and $\widehat{\phi}_t$, and (ii) that $\phi_x$ and $\phi_t$ decay to $\tau$ at the first interior points of their supports in real space (as in \eqref{realdecay}). For (i) we utilize that test functions $\phi_{a,p}\in \CalS$ defined in \eqref{testfcn} of the form 
\[\phi_{a,p}(s) = C\left(1-\left(\frac{s}{a}\right)^2\right)^p_+\]
are well-approximated by Gaussians. Indeed, letting $C$ be such that $\nrm{\phi_{a,p}}_1=1$ and setting $\sigma := a/\sqrt{2p+3}$, then $\phi_{a,p}$ matches the first three moments of the Gaussian
\[\rho_\sigma(s):= \frac{1}{\sqrt{2\pi \sigma^2}}e^{-s^2/2\sigma^2},\] 
which provides a bound on the error in the Fourier transforms $\widehat{\phi}_{a,p}$ and $\widehat{\rho}_\sigma$ for small frequencies $\xi$ in terms of their 4th moments\footnote{This also shows that with $\sigma = a/\sqrt{2p+3}$, if we take $a = \sqrt{2p}$ then we get pointwise convergence $\phi_{a,p}\to \rho_1$ as $p\to \infty$.}:
\[|\widehat{\phi}_{a,p}(\xi)-\widehat{\rho}_\sigma(\xi)| \leq |\xi|^4 \left(\frac{a^4}{2}\left[\frac{p+3/2}{(4p^2+12p+9)(4p^2+16p+15)}\right] + o(1)\right) = \CalO(|\xi|^4a^4p^{-3}).\]
For small $\xi$ and $a$ and large $p$, it suffices to use $\widehat{\rho}_\sigma(\xi) = \rho_{1/\sigma}(\xi)$ as a proxy for $\widehat{\phi}_{a,p}$.  

To enforce decay of $\phi_x$ in Fourier space (and similarly for $\phi_t$) we specify that $k^*_x$ is $\widehat{\tau}$ standard deviations into the tail of $\rho_{1/\sigma}(\xi)$, where $\sigma = a/\sqrt{2p+3}$. Recalling that $a = m_x\Delta x = m_x(L/N_1)$ where $L$ is the length of the spatial domain and $N_1$ is the number of points in $x$, this provides a relation between the degree $p_x$ and the discrete support hyperparameter $m_x$:
\[\frac{2\pi}{L}k^*_x = \frac{\widehat{\tau}}{\sigma} = \widehat{\tau}\frac{\sqrt{2p_x+3}}{a}= \widehat{\tau}\frac{\sqrt{2p_x+3}}{m_x(L/N_1)}\]
\[\symb{0.4}{\implies}  2\pi k^*_xm_x = \widehat{\tau}N_1\sqrt{2p_x+3}.\]
Enforcing decay in real space then provides a second condition between $p_x$ and $m_x$: 
\[\left(1-(1-1/m_x)^2\right)^{p_x} = \tau.\]
Combining these two constraints we get that $m_x$ will be a root of 
\[F(m) := F(m;\ k_x,N_1,\widehat{\tau},\tau) := \log\left(\frac{2m-1}{m^2}\right)\left(4\pi^2{k_x^*}^2m^2-3N_1^2\widehat{\tau}^2\right)-2N_1^2\widehat{\tau}^2\log(\tau).\]
Provided $N_1>4$, $0<\tau<1$ and $\frac{2\pi}{\sqrt{3}}\left(\frac{k^*_x}{N_1/2}\right)\leq \widehat{\tau}\leq \frac{\pi}{\sqrt{3}}k_x^*$, then $F(m)$ has a unique root $m_x\geq 2$ in the nonempty interval 
\[\Bigg[\frac{\sqrt{3}}{\pi}\left(\frac{N_1/2}{k_x^*}\right)\widehat{\tau},\ \ \, \frac{\sqrt{3}}{\pi}\left(\frac{N_1/2}{k_x^*}\right)\widehat{\tau}\sqrt{1-(8/\sqrt{3})\log(\tau)}\Bigg]\] 
on which $F$ is monotonically decreasing and changes sign. After finding $m_x$ we can solve for $p_x$ using either constraint. Figure \ref{visKS} illustrates this process for computing the column of $\Gbf$ corresponding to the Burgers-type nonlinearity $\partial_x(u^2)$ using the same KS dataset as in Figure \ref{visKS} with $50\%$ noise. A one-dimensional slice in $x$ is taken at fixed time $t = 99$ and compared with the underlying clean data, showing that the convolution successfully filters out the noise-dominated modes (despite the fact that the noisy, nonlinearly-transformed data $(\Ubf)^2$ has a substantial bias compared with the clean data $(\Ubf^\star)^2$).
%Ultimately this trades the search for $m$ and $p$ values to specify $\phi$ into a search for $J$, the number of standard deviations into the tail of $\phi_{a,p}$ that we wish to place the noise-dominated region of $\CalF(\Ubf)$ (note that $\tau$, the decay tolerance, was already a hyperparameter in the WSINDy\_PDE and it doesn't seem to have a large effect so long as $\tau\in [10^{-8}, 10^{-16}]$). Is this any better? Sure: our choice in $J$ exactly translates to how much weight we are giving to the noise-dominated part of the spectrum. In particular, if a model has a $u^2$ nonlinearity, then $\widehat{u^2}(\xi) = \widehat{u}*\widehat{u}(\xi) \approx \widehat{u}(\xi/2)$, so perhaps we don't want to pick $J$ to be too large. Consequently, for nonlinear Schr\"odinger's and reaction-diffusion, $J = 2/3$ leads to good recovery, seemingly because of cubic nonlinearities, while for Burgers-type nonlinearities $J=2$ suffices, as $\widehat{u^2}\cdot \widehat{\phi_x}$ ``spreads out'' the spectra of both $u$ and $\phi$, and so if we place the noise-dominated spectrum of $\Ubf$ two standard deviations $(J=2)$ into the tail of $\widehat{\phi}$, the the noise-dominated modes of ${\Ubf}^2$ also land in the tail of $\widehat{\phi_x}$. 

\section{Numerical Simulation Methods}\label{app:numpde}

We now review the numerical methods used to simulate noise free data sets for each of the PDEs in Table \ref{pdetable} (note that dimensions of the datasets are given in Table \ref{hptable}). With the exception of the Navier-Stokes equations, which was simulated using the immersed boundary projection method in C++ \cite{TaiCol-06}, all computations were performed in MATLAB 2019b. 

\subsection{Inviscid Burgers}
\begin{equation}\label{burgs}
\partial_t u  = -\frac{1}{2} \partial_x (u^2)
\end{equation}
We take for exact data the shock-forming solution
\begin{equation}\label{burgers}
u(x,t) = \begin{dcases} A, & t\geq \max\left\{\frac{1}{A}x+\frac{1}{\alpha},\ \frac{2}{A}x+\frac{1}{\alpha}\right\} \\ -\frac{\alpha x}{1-\alpha t}, & A\left(t - \frac{1}{\alpha}\right)<x\leq 0 \\ 0,& \text{otherwise}\end{dcases}.
\end{equation}
which becomes discontinuous at $t= \alpha^{-1}$ with a shock travelling along $x = \frac{A}{2}\left(t-\frac{1}{\alpha}\right)$ (see Figure \ref{burgerssols}). We choose $\alpha = 0.5$ and an extreme value of $A = 1000$ to demonstrate that WSINDy still has excellent performance for large amplitude data. The noise-free data consists of $\eqref{burgers}$ evaluated at the points $(x_i,t_j) = (-4000+i\Delta x, j\Delta t)$ with $\Delta x = 31.25$ and $\Delta t = 0.0157$ for $1\leq i,j\leq  256$. 

\subsection{Korteweg-de Vries}\label{sec:KdV}
\begin{equation}\label{KdV}
\partial_t u  = -\frac{1}{2} \partial_x (u^2) -\partial_{xxx}u
\end{equation}
A solution is obtained for $(x,t) \in [-\pi,\pi]\times[0, 0.006]$ with periodic boundary conditions using ETDRK4 timestepping and Fourier-spectral differentiation \cite{kassam2005fourth} with $N_1 = 400$ points in space and $N_2 = 2400$ points in time. We subsample $25\%$ of the timepoints for system identification and keep all of the spatial points for a final resolution of $\Delta x = 0.0157$, $\Delta t = 10^{-5}$. For initial conditions we use the two-soliton solution
\[u(x,0) = 3A^2\text{sech}(0.5(A(x+2)))^2 + 3B^2\text{sech}(0.5(B(x+1)))^2, \qquad A=25, B=16.\] 
\subsection{Kuramoto-Sivashinsky}\label{sec:KS}
\begin{equation}\label{KS}
\partial_t u  = -\frac{1}{2} \partial_x (u^2) -\partial_{xx}u -\partial_{xxxx} u.
\end{equation}
A solution is obtained for $(x,t) \in [0,32\pi]\times[0, 150]$ with periodic boundary conditions using ETDRK4 timestepping and Fourier-spectral differentiation \cite{kassam2005fourth} with $N_1 = 256$ points in space and $N_2 = 1500$ points in time. For system identification we subsample $20\%$ of the time points for a final resolution of $\Delta x = 0.393$ and $\Delta t = 0.5$. For initial conditions we use
\[u(x,0) = \cos(x/16)(1+\sin(x/16)).\] 

\subsection{Nonlinear Schr\"odinger}\label{sec:NLS}
\begin{equation}\label{NLS}
w_t = -\frac{i}{2}\partial_{xx}w + |w|^2w
\end{equation}
For the nonlinear Schr\"odinger equation (NLS) we reuse the same dataset from \cite{rudy2017data}, containing $N_1=512$ points in space and $N_2 = 502$ timepoints, although we subsample $50\%$ of the spatial points and $50\%$ of the time points for a final resolution of $\Delta x = 0.039$, $\Delta t = 0.0125$. For system identification, we break the data into real and imaginary parts $(w = u+iv)$ and recover the system
\begin{equation}\label{NLSsep}
\begin{cases}\partial_t u= \frac{1}{2}\partial_{xx}v +u^2v+v^3 \\ \partial_t v = -\frac{1}{2}\partial_{xx} u-uv^2-u^3.\end{cases}
\end{equation}
\subsection{Sine-Gordon}\label{sec:SG}
\begin{equation}\label{SG}
\partial_{tt} u  = \vphantom{\frac11} \partial_{xx}u+\partial_{yy}u - \sin(u)
\end{equation}
A numerical solution is obtained using a pseudospectral method on the spatial domain $[-\pi,\pi]\times [-1,1]$ with 64 equally-spaced points in $x$ and 64 Legendre nodes in $y$. Periodic boundary conditions are enforced in $x$ and homogeneous Dirichlet boundaries in $y$. Geometrically, waves can be thought of as propagating on a right cylindrical sheet with clamped ends. Leapfrog time-stepping is used to generate the solution until $T=5$ with $\Delta t = 6$e$-5$. We then subsample $0.25\%$ of the timepoints and interpolate onto a uniform grid in space with $N_1=403$ points in $x$ and $N_2 = 129$ points in $y$. The final resolution is $\Delta x = 0.0156$, $\Delta t = 0.025$. We arbitrarily use Gaussian data for the initial wave disturbance:
\[u(x,y,0) = 2\pi\exp(-8(x-0.5)^2 - 8y^2).\]

\subsection{Reaction-Diffusion}\label{sec:RD}

\begin{equation}\label{RD}
\begin{cases} \partial_t u  = \frac{1}{10} \partial_{xx}u+\frac{1}{10} \partial_{yy}u -uv^2-u^3+v^3+u^2v+u \\ \partial_t v  = \frac{1}{10} \partial_{xx}v+\frac{1}{10} \partial_{yy}v+v-uv^2-u^3-v^3-u^2v \end{cases}
\end{equation}

The system \eqref{RD} is simulated over a doubly-periodic domain $(x,y)\in[-10,10]\times[-10,10]$ with $t\in [0,10]$ using Fourier-spectral differentiation in space and method-of-lines time integration via MATLAB's $\texttt{ode45}$ with default tolerance. The computational domain has dimensions $N_1=N_2=256$ and $N_3=201$, for a final resolution of $\Delta x = 0.078$, $\Delta t = 0.0498$. For initial conditions we use the spiral data 
\[\begin{dcases} u(x,y,0)=\text{tanh}(\sqrt{x^2+y^2})\cos\left(\theta(x+iy)-\pi\sqrt{x^2+y^2}\right)\\
v(x,y,0)=\text{tanh}(\sqrt{x^2+y^2})\sin\left(\theta(x+iy)-\pi\sqrt{x^2+y^2}\right),\end{dcases}\]
where $\theta(z)$ is the principle angle of $z\in \Cbb$. Note that this is an unstable spiral which breaks apart over time but still settles into a limit cycle. 

\begin{figure}
\begin{tabular}{cc}
		\includegraphics[trim={10 0 25 20},clip,width=0.45\textwidth]{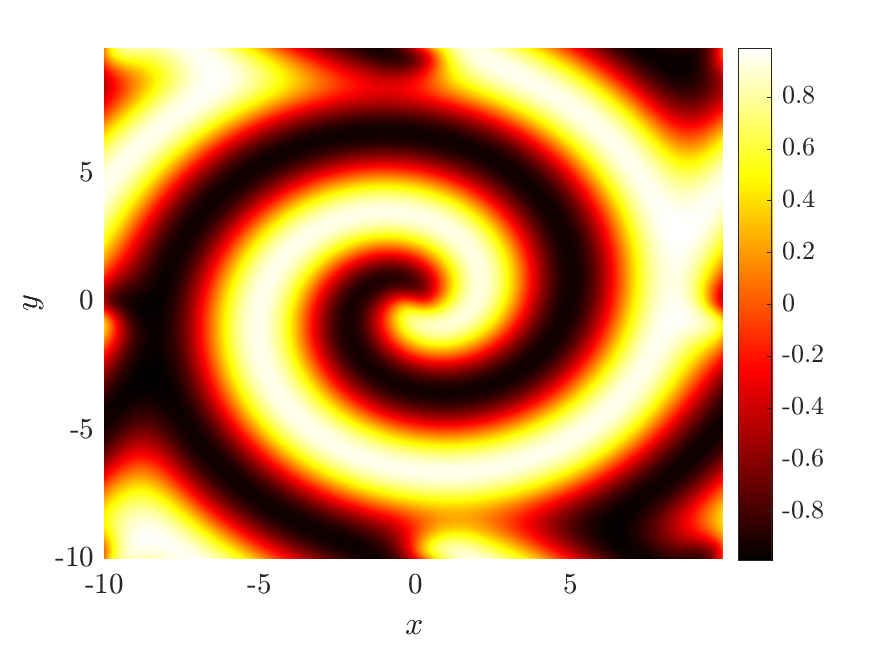} &
		\includegraphics[trim={10 0 25 20},clip,width=0.45\textwidth]{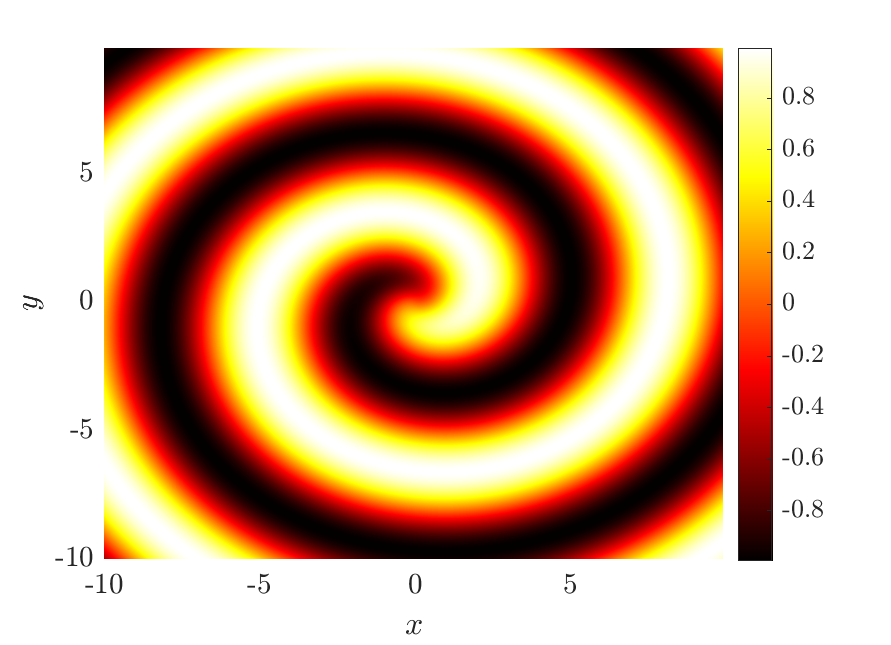}  
\end{tabular}
\caption{Comparison between the full reaction-diffusion model \eqref{RD} (left) and the pure-oscillatory reduced model \eqref{RDredux} (right) at the final time $T=10$ with both models simulated from the same initial conditions leading to a spiral wave (only the $v$ component is shown, results for $u$ are similar). The reduced model provides a good approximation away from the boundaries.}
\label{rdcompare}
\end{figure}

Using the traditional (stable) spiral wave data \cite{rudy2017data} (differing only from the dataset used here in that the term $\pi\sqrt{x^2+y^2}$ in the initial conditions above is replaced by $\sqrt{x^2+y^2}$) we noticed an interesting behavior in that for high noise the resulting model was purely oscillatory. In other words, the stable spiral limit cycle  happens to be well-approximated by the pure-oscillatory model
\begin{equation}\label{RDredux}
\partial_t\ubf = \alpha\begin{pmatrix}0 & 1 \\ -1 & 0 \end{pmatrix}\ubf
\end{equation}
with $\alpha \approx 0.91496$. A comparison between this purely oscillatory reduced model and the full model simulated from the same initial conditions is shown in Figure \ref{rdcompare}. For $\sigma_{NR}\leq 0.1$ WSINDy applied to the stable spiral dataset returns the full model, while for $\sigma_{NR}> 0.1$ the oscillatory reduced model is detected. This suggests that although the stable spiral wave is a hallmark of the $\lambda$-$\omega$ reaction-diffusion system, from the perspective of data-driven model selection it is not an ideal candidate for identification of the full model.

\subsection{Navier-Stokes}\label{sec:NS}

\begin{equation}\label{NS}
\partial_t \omega = -\partial_x(\omega u)-\partial_y(\omega u) +\frac{1}{100} \partial_{xx}\omega+\frac{1}{100} \partial_{yy}\omega
\end{equation}
A solution is obtained on a spatial grid $(x,y) \subset [-1,8]\times[-2,2]$ with a ``cylinder'' of diameter $1$ located at $(0,0)$. The immersed boundary projection method \cite{TaiCol-06} with 3rd-order Runge-Kutta timestepping is used to simulate the flow at spatial and temporal resolutions $\Delta x=\Delta t = 0.02$ for 2000 timesteps following the onset of the vortex shedding limit cycle. The dataset $(\Ubf,\Vbf,\Wbf)$ contains the velocity components as well as the vorticity for points away from the cylinder and boundaries in the rectangle $(x,y)\in [1,7.5]\times[-1.5,1.5]$. We subsample $10\%$ of the data in time for a final resolution of $\Delta x = 0.02$ and $\Delta t = 0.2$. 

\end{document}